# Gamifying optimization: a Wasserstein distance-based analysis of human search

*Antonio Candelieri [1], Andrea Ponti [2], Francesco Archetti [2]*

[1] University of Milano-Bicocca, Department of Economics, Management and Statistics, Milan, Italy
[2] University of Milano-Bicocca, Department of Computer Science, Systems and Communication, Milan, Italy

Corresponding author: antonio.candelieri@unimib.it
ORCID: Antonio Candelieri (0000-0003-1431-576X), Francesco Archetti (0000-0003-1131-3830), AndreaPonti (0000-0003-4187-4209).

**Abstract**. The main objective of this paper is to outline a theoretical framework to characterise humans' decision-making strategies under uncertainty, in particular active learning in a black-box optimization task and trading-off between information gathering (exploration) and reward seeking (exploitation). Humans' decisions making according to these two objectives can be modelled in terms of Pareto rationality. If a decision set contains a Pareto efficient (dominant) strategy, a rational decision maker should always select the dominant strategy over its dominated alternatives. A distance from the Pareto frontier determines whether a choice is (Pareto) rational. To collect data about humans' strategies we have used a gaming application that shows the game field, with previous decisions and observations, as well as the score obtained.
The key element in this paper is the representation of behavioural patterns of human learners as a discrete probability distribution, specifically a histogram. In this representation the similarity between users' behaviours can be captured by a distance between their associated histograms. This maps the problem of the characterization of humans' behaviour into a space whose elements are probability distributions structured by a distance between histograms, namely the transport-based Wasserstein distance. The distributional analysis gives new insights about human search strategies and their deviations from Pareto rationality. Since the uncertainty is one of the two objectives defining the Pareto frontier, the analysis has been performed for three different uncertainty quantification measures to identify which better explains the Pareto compliant behavioural patterns. Beside the analysis of individual patterns Wasserstein has also enabled a global analysis computing the barycenters and Wasserstein k-means clustering. A further analysis has been performed by a decision tree to relate non-Paretian behaviour, characterized by "exasperated exploitation", to the dynamics of the evolution of the reward seeking process.

**Acknowledgments.** We greatly acknowledge the DEMS Data Science Lab, Department of Economics Management and Statistics (DEMS), for supporting this work by providing computational resources.

# 1. Introduction

## 1.1 Motivation

Human activity at all levels and time scales of decision-making under uncertainty requires balancing e*xploitation* – meaning the use of the *knowledge* collected so far to maximize immediate reward – and *exploration* – meaning investing resources to acquire more knowledge to update one's beliefs. This balance is usually called the exploration-exploitation dilemma: decisions allowing for increasing knowledge do not necessarily lead to the greatest immediate reward (Wilson et al., 2020a; Wilson et al., 2014).
The trade-off between explorative and exploitative behaviours characterizes many disciplines (Berger-Tal et al., 2014) and has originated a multidisciplinary framework that applies to humans, animals, and organizations. The analysis of the strategies implemented by humans in dealing with uncertainty has been an actively researched topic (Schulz et al., 2015; Gershman, 2018; Schulz & Gershman, 2019). A key observation, motivating this line of research, is also the awareness that *human learners* are amazingly fast and effective at adapting to unfamiliar environments and incorporating upcoming knowledge: this is an intriguing behaviour for cognitive sciences as well as an important challenge for Machine Learning. The reference task considered in this paper is the optimization problem:

$$x^* = \underset{x\in\Omega\subset\Re^d}{\mathrm{argmax}} f(x) \qquad (1)$$

with $f(x)$ black box, meaning that its analytical form is not given, no derivatives are available, and the value of $f(x)$ can be only known pointwise through expensive and noisy evaluations. Finally, $\Omega$ denotes the *search space*, usually box bounded.

We consider sequential optimization to solve (1). At a generic iteration $n$, the player/agent/algorithm chooses a location $x^{(n)}$ to query and observe/collect the associated function value, possibly perturbed by noise, that is $y^{(n)} = f(x^{(n)}) + \varepsilon$. The goal is to get close to the optimizer $x^*$ within a limited number, $N$, of trials. The choice is performed according to a so-called acquisition (or infill) function which, at each iteration, trades-off between exploration and exploitation. A related goal is to maximize the Average Cumulative Reward (ACR) over the $N$ trials, that is $ACR^{(N)} = \frac{1}{N}\sum_{i=1}^{N} y^{(i)}$.

Recently, the Bayesian Optimization framework (BO) (Frazier, 2018; Archetti & Candelieri, 2019) has become one of the most efficient method for solving (1). This relates to the ongoing discussion in cognitive science as to whether also humans' strategies are sample efficient: (Borji & Itti, 2013; Candelieri et al., 2020) have been arguing, based on empirical evidence, that strategies adopted by humans in solving global optimization problems have a much stronger association with BO than to other optimization algorithms. Gaussian Process (GP) modelling and Bayesian learning, first proposed in (Kruschke, 2008; Griffiths et al., 2008) have emerged as central paradigms in modelling human learning, where the GP model is used to approximate the outcome of the next decision conditioned on previous decisions and observed outcomes. Fitting a GP requires to choose, a priori, a *kernel* as covariance function which determines the predictive uncertainty; different kernels are available, each one implying a different characterization for the approximation of $f(x)$. As already stated in (Wilson et al., 2015), it was demonstrated that *"GPs with standard kernels struggle on function extrapolation problems that are trivial for human learners"*. A kernel learning framework is proposed to reverse engineer the inductive biases of human learners across a set of behavioural experiments, gaining psychological insights and extrapolating in humanlike ways that go beyond traditional kernels. Moreover, (Gershman, 2019) remarked that different quantifications of the uncertainty – as discussed later in the paper – are a key concept also in theories of cognition and emotion.

The key element in this paper is the representation of behavioural patterns of human learners, over different tasks, as a discrete probability distribution, specifically a histogram. A behavioural signature is associated to everyone, described by an estimated probability measure like a histogram expressed by a sequence of disjoint intervals with weights. In this sense our approach is related to Symbolic Data Analysis (SDA) (Bock & Diday, 2012). In this representation the similarity between users' behaviours can be captured by a distance between their associated histograms. This maps the problem of the characterization of humans' behaviour into a space whose elements are histograms: this space is structured by a distance between histograms, namely the Wasserstein distance.

In this paper we focus on the Wasserstein (WST) distance, which is based on a different approach, optimal transport, which is a field of mathematics which studies the geometry of probability spaces and provides a principled framework to compare and align probability distributions. The Wasserstein distance can be traced back to the works of Gaspard Monge (1781) and Lev Kantorovich (1942). Specifically, we consider one instance of WST, called Earth Mover Distance (EMD), which is a natural and intuitive distance between discrete probability distributions and in particular histograms. Let's think of histograms as piles of sand sitting on the ground (underlying domain). To measure the difference between distributions we measure how far the grains of sand have to be moved so that the two distributions exactly coincide.

Defining a distance between distributions require a notion of distance between points in the underlying domain which is called the ground distance. EMD is the minimum ground distance travelled weighted by the amount of sand moved (called flow). If the ground distance is a metric and distributions have the same mass, which is the case of probability distribution functions, EMD is a metric as well. Also, multi-dimensional discrete distributions can be defined as histograms by partitioning the underlying domain into bins with a weight associated to each bin.

WST has evolved into a very rich mathematical structure whose complexity and flexibility are analysed in a landmark volume (Villani, 2009).

Computing Wasserstein distances requires in general the solution of a constrained linear optimization problem which has a very large number of variables and constraints when the support of the probability distributions is multidimensional. It is shown to be equivalent to a min-flow problem. Recently, many specialized approaches have drastically reduced the computational hurdles (Peyré & Cuturi, 2019).

The main advantage WST is that it is a cross binning distance, and it is not affected by different binning schemes. Moreover, WST matches naturally the perceptual notion of nearness and similarity. This is not the case of KL and $\chi$-square distances that account only for the correspondence between bins of the same index and do not use information across bins or distributions with different binning schemes, that is different support. Two important elements of the WST theory are the barycenter and the Wasserstein clustering. The Wasserstein barycenter offers a useful characterization of a set of distributions. Also, a standard clustering method like k-means can be generalized to WST spaces. Beside the analysis of individual patterns, these tools have also enabled a global analysis computing the WST barycenters and k-mean WST clustering, and will be used, as shown in Sect. 6, to characterize and classify the behavioural patterns.

The analysis of the Pareto frontier in the space of the GP's predictive mean and standard deviation offers a set of Pareto-efficient decisions which can be significantly wider than those selected through "traditional" acquisition functions. According to Pareto-based rationality model, if a decision set contains a Pareto efficient (dominant) strategy, a rational decision maker should always select the dominant strategy over its dominated alternatives. Still, according to a famous Schumpeter quotation (Schumpeter, 1954) traditional decision making under risk *"has a much better claim to being called a logic of choice than a psychology of value"* and indeed deviations from Pareto rational behaviour have been documented in domains like economics, business, but also Reinforcement Learning.

The analysis of violations of dominance in decision-making has become mainstream economics under the name of behavioural economics and prospect theory (Kahneman, 2011): rather than being labelled "irrational", non-Pareto compliant behaviour is just not well described by the rational-agent model.

Since the uncertainty is one of the two objectives defining the Pareto frontier, the analysis has been performed for three different uncertainty quantification measures with the aim to identify the one resulting more compliant with the Pareto rationality model. Beside the analysis of individual patterns Wasserstein has also enable a global analysis computing the WST barycenters and WST k-mean clustering.

The key result is an analytical framework to characterize how deviations from "rationality" depend on individual features represented in the histogram and uncertainty quantifications.

A hypothesis that we also wanted to test is that a sizable number on non-Paretian instances can be associated to "exasperated exploration" related to the evolution of the reward seeking process expressed by the dynamics of Average Cumulated Reward (ACR): this has been analysed by a decision tree classifier in which the "class label" (Paretian vs non-Paretian) is predicted.

## 1.2 Contributions of this paper

The behavioural data have been analysed also according to a Pareto model: human search data have shown significant deviations from "rationality" depending on individual features, represented as histograms, and uncertainty quantifications. The key contribution of the paper is the proposal of a distributional analysis of human search pattern based on the Wasserstein distance. This distributional analysis has been conducted at the individual level and an aggregate level computing barycenters and performing clustering in the Wasserstein space. It is also interesting to remark that while most of the previous works addressed how people assess the information value of possible queries, in this paper we rather address the issue of the perception of probabilistic uncertainty itself. Note that BO algorithm is not actually executed: sequences of points are generated by humans and compared with optimal Pareto fronts generated analytically.

An important contribution has been the development of a software environment for gathering data about human behaviour and analysing them, whose use can be helpful, beyond the specific case, to analyse human strategies in learning problems.

Another contribution has been the development of a decision tree classifier used to characterize deviations from Pareto "rationality" and switches towards "exasperated" exploration depending on the subject, the specific task at hand, and the evolution of the optimization process, measured as the reward collected over players' decisions.

The computational results and their analysis allow to formulate at least a tentative answer to the following research questions:
- Do humans always make "rational" choices (i.e., Pareto optimal decisions between the improvement expected and uncertainty) or, in some cases, they "exasperate" exploration?
- Do different uncertainty quantification measures lead to different classifications of humans' decisions? And which uncertainty quantification measure make humans "more rational"?
- Does the distributional representation provide an efficient signature of the subject/task?
- Does the WST distance capture the difference between the behaviour of 2 subjects on a given task?
- What's the average behaviour on a given task of all the subjects?
- What is an index of the difficulty of a task averaged over all subjects?
- What is the relative impact on the Pareto compliance of a human of the distance and the kernel used?

**1.3 Related works**

In Sect. 1.1 we have briefly introduced the issue of uncertainty quantification in humans and its relationship with learning and optimization and new analytical tools to characterize humans' behaviour. Here we provide a more specific analysis of the prior work and significant recent results carried out for Wasserstein.

*1.3.1 Wasserstein*
The literature on WST is now immense: general references have been already given in the introduction, other will be given in Sec. 4. Here we limit to very specific references mainly focused on WST clustering. As this is a very demanding computational task, this paper addresses this issue by proposing more efficient methods. An early contribution is from (Applegate et al., 2011), who proposes an EMD-based clustering to analyze mobility usage patterns. EMD is shown to allow for representing inherent relationships in the WST space and to cluster meaningfully also sparse signatures. (Cabanes et al., 2021) uses a dimensionality reduction by Self-organizing Maps (SOM) learning and then cluster data within a WST space. (Puccetti et al., 2020) proposes an Iterative Swapping Algorithm (ISA) which is shown to have a quadratic complexity in the barycenter computation. (Ye et al., 2017) proposes an approach based on the Alternating Direction Method of Multipliers (ADMM) for WST clustering of distribution with sparse support. (Verdinelli & Wasserman, 2019) introduces a hybrid distance based on gaussian approximations.

*1.3.2. Cognitive Sciences*
An early contribution (Cohen et al., 2007) analyses how humans manage the trade-off between exploration and exploitation in non-stationary environments. Successively, (Wilson et al., 2014) demonstrates that humans use both *random* and *directed exploration*. (Gershman & Uchida, 2019) show how directed exploration in humans amounts to adding an *"uncertainty bonus"* to estimated reward values and how this brings to the *Upper Confidence Bound* (UCB) acquisition function in Multi Armed Bandits (Auer et al., 2002) and BO (Srinivas et al., 2012). The same approach is elaborated in (Schulz & Gershman, 2019), who distinguish between *irreducible uncertainty*, related to the reward stochasticity, and *uncertainty*, which can be reduced through information gathering. In the former the decision strategy is *random search* while for the latter is *directed exploration* which attaches an uncertainty bonus to each decision value. This distinction mirrors the one in Machine Learning between *aleatoric* uncertainty – due to the stochastic variability inherent in querying $f(x)$ – and *epistemic* uncertainty – due to the lack of knowledge about the actual structure of $f(x)$ – which can be reduced by collecting more information. (Friston et al., 2014) analyses how entropy and expected utility account, respectively, for exploratory and exploitative behaviours, "*relating it to the discrepancy between observed and expected reward, known as the reward prediction error (RPE), which serves as a learning signal for updating reward expectations. On the other hand, dopamine also appears to participate in various probabilistic computations, including the encoding of uncertainty and the control of uncertainty-guided exploration*" (Gershman & Uchida, 2019).

*1.3.3. Bayesian Optimization*
In the BO research community, recent papers proposed to analyze the exploration-exploitation dilemma as a bi-objective optimization problem: minimizing the predictive mean (associated to exploitation) while maximizing uncertainty, typically the predictive standard deviation as in UCB (associated to exploration). (Žilinskas & Calvin, 2019). This mean-variance framework has been also considered in (Iwazaki et al., 2021), for multi-task, multi-objective, and constrained optimization scenarios. (De Ath et al., 2019; De Ath et al., 2020) show that taking a decision by randomly sampling from the Pareto frontier can outperform other acquisition functions. The main motivation is that the Pareto frontier offers a set of Pareto-efficient decisions wider than that allowed by "traditional" acquisition functions.

*1.3.4. Economics*
The issue of deviations from Pareto optimality has become a central topic in behavioural economics from the seminal work in (Tversky & Kahneman, 1989) to (Kourouxous & Bauer, 2019) which identifies the most common causes for violations of dominance, namely *framing* (i.e., presentation of a decision problem), *reference points* (i.e., a form of prior expectation), *bounded rationality* and *emotional responses*. An entirely different approach is suggested in (Peters, 2019) where the concept of ergodicity from statistical mechanics is proposed to model non Paretian behaviour.

A recent important contribution is (Sandholtz, 2020) which tackles the problem to infer, given the observed search path generated by a human subject in the execution of a black box optimization task, the unknown acquisition function underlying the sequence. For the solution of this problem, referred to as Inverse Bayesian Optimization (IBO), a probabilistic framework for the non-parametric Bayesian inference of the acquisition function is proposed, performed on a set of possible acquisition functions.

**1.4 Outline of the paper**

Sect. 2 introduces the definitions of Gaussian Process regression and different uncertainty quantifications. Sect. 3 introduces the Wasserstein distance, both the basic notions and the computational issues. Sect. 4 develops a framework for the application of the Pareto analysis to the specific problem considering different kernels for the GP regression models as well as three different uncertainty quantification measures. Sect. 5 introduces the experimental framework used for data collection, that is the decisions taken by the human players according to their personal search strategies, and the proposed analytical framework. Sect. 6 describes the relevant results obtained by the application of the analytical framework. Finally, Sect. 7 outlines the conclusions of this study and the perspective of future works.

In order not to break down the chain of arguments, most of the background material and the description of software and data resources are given in the supplementary material.

## 2. Materials and methods

### 2.1 Gaussian Process regression

A GP is a *random distribution over functions* $f: \Omega \subset \Re^d \to \Re$ denoted with $f(x) \sim GP(\mu(x), k(x,x'))$ where $\mu(x) = \mathbb{E}(f(x)): \Omega \to \Re$ is the mean function of the GP and $k(x,x'): \Omega \times \Omega \to \Re$ is the *kernel* or *covariance function*. One way to interpret a GP is as a collection of correlated random variables, any finite number of which have a joint Gaussian distribution, so $f(x)$ can be considered as a sample drawn from a multivariate normal distribution. In Machine Learning, GP modelling is largely used for both classification and regression tasks (Williams & Rasmussen, 2006; Gramacy, 2020), providing probabilistic predictions by conditioning $\mu(x)$ and $\sigma^2(x)$ on a set of available data/observations.

Let denote with $X_{1:n} = \{x^{(i)}\}_{i=1,\ldots,n}$ a set of $n$ locations in $\Omega \subset \Re^d$ and with $y_{1:n} = \{f(x^{(i)}) + \varepsilon\}_{i=1,\ldots,n}$ the associated function values, possibly noisy with $\varepsilon$ a zero-mean Gaussian noise $\varepsilon \sim \mathcal{N}(0, \lambda^2)$. Then $\mu(x)$ and $\sigma^2(x)$ are the GP's posterior predictive mean and standard deviation, conditioned on $X_{1:n}$ and $y_{1:n}$ according to the following equations:

$$\mu(x) = k(x, X_{1:n}) \, [K + \lambda^2 I]^{-1} \, y_{1:n} \qquad (2)$$

$$\sigma^2(x) = k(x, x) - k(x, X_{1:n}) \, [K + \lambda^2 I]^{-1} \, k(X_{1:n}, x) \qquad (3)$$

where $k(x, X_{1:n}) = \{k(x, x^{(i)})\}_{i=1,\ldots,n}$ and $K \in \Re^{n \times n}$ with entries $K_{ij} = k(x^{(i)}, x^{(j)})$.

The choice of the kernel establishes prior assumptions over the structural properties of the underlying (aka latent) function $f(x)$, specifically its smoothness. However, almost every kernel has its own hyperparameters to tune – usually via Maximum Log-likelihood Estimation (MLE) or Maximum A Posteriori (MAP) – for reducing the potential mismatches between prior smoothness assumptions and the observed data. Common kernels for GP regression – considered in this paper – are:

- Squared Exponential: $k_{SE}(x, x') = e^{-\frac{\|x - x'\|^2}{2\ell^2}}$
- Exponential: $k_{EXP}(x, x') = e^{-\frac{\|x - x'\|}{\ell}}$
- Power-exponential: $k_{PE}(x, x') = e^{-\frac{\|x - x'\|^p}{\ell^p}}$
- Matérn3/2: $k_{M3/2}(x, x') = \left(1 + \frac{\sqrt{3}\,\|x - x'\|}{\ell}\right) e^{-\frac{\sqrt{3}\,\|x - x'\|}{\ell}}$
- Matérn5/2: $k_{M5/2}(x, x') = \left[1 + \frac{\sqrt{5}\,\|x - x'\|}{\ell} + \frac{5}{3}\left(\frac{\|x - x'\|}{\ell}\right)^2\right] e^{-\frac{\sqrt{5}\,\|x - x'\|}{\ell}}$

The main well-known disadvantage of GP modelling is its cubic complexity due to the inversion of the matrix $[K + \lambda^2 I]$.

## 2.2 Uncertainty quantification and active learning

From a modelling perspective, uncertainty can be split in aleatoric and epistemic (Der Kiureghian & Ditlevsen, 2009; Kendall & Gal, 2017), where the aleatoric uncertainty is randomness proper in the valuation of $f(x)$ (usually named "noise") and cannot be reduced, while the epistemic uncertainty depends on the model and can be reduced by collecting more data.

From an informal – yet more intuitive – point of view, uncertainty about a decision is the amount of lack of knowledge about it, increasing with the "distance" from decisions already performed and where "distance" can be any suitable metric to compare two decisions. When decisions are locations in a search space, as in this paper, any spatial distance can be considered: an example of this uncertainty quantification has been recently proposed in (Bemporad, 2020) which uses Radial Basis Functions (RBF) as surrogate model and an inverse distance weighting such that the proposed distance is zero at sampled points and grows in between

Finally, when humans' decisions are analysed, there is still another relevant lack in mathematical methods for uncertainty quantification as recently demonstrated in (Schulz et al., 2019) and (Bertram et al., 2020), which have investigated the role of emotion in judgment, risk assessment, and decision making under uncertainty and the different kinds of entropy which can be used to quantify uncertainty and shown that emotional states are significantly connected with subjective uncertainty estimation.

From the viewpoint of Machine Learning, uncertainty quantification plays a pivotal role in reduction of errors during learning, optimization and decision making. In (Abdar et al., 2021) a wide survey of different uncertainty quantification methods is provided, considering many application fields.

In decision making, uncertainty is usually associated to exploration: when the uncertainty is "larger" than the possible estimated improvement, then it could be more profitable to adopt an explorative behaviour and acquire more knowledge about $f(x)$.

Acquisition functions are analysed in appendix Al. In what follows only few basic definitions.

Global optimization methods differ one from another in how they generate the next decision (i.e., location) $x^{(n+1)}$. To do this, BO fits a GP according to (2-3) and where $X_{1:n} = \{x^{(i)}\}_{i=1,\ldots,n}$ and $y_{1:n} = \{y^{(i)}\}_{i=1,\ldots,n}$ are the two sequences of, respectively, decisions made and associated observed outcomes. Then, an acquisition function, combining GP's $\mu(x)$ and $\sigma(x)$, is optimized to obtain $x^{(n+1)}$, while dealing with the exploration-exploitation trade-off. GP-UCB, it is also classified as an optimistic policy, because it chooses $x^{(i+1)}$ depending on the most optimistic value for $f(x)$ under the current GP. From a cognitive point of view, (Wu et al., 2018) analysed the human search strategy, under a limited number of trials, concluding that GP-UCB offers the best option for modelling the exploitation-exploration trade-off adopted by the humans.

Let $\mathcal{K}$ denotes the set of kernels to choose as GP's prior. In this study $\mathcal{K} = \{k_{SE}, k_{EXP}, k_{PE}, k_{M3/2}, k_{M5/2}\}$.
Let $\zeta(x)$ denotes the improvement expected by querying the objective function at location $x$, depending on the GPs' posterior (i.e., one GP for each kernel in $\mathcal{K}$). Formally, $\zeta(x) = \mu(x) - y^+$, where $y^+ = \max_{i=1,\ldots,n} \{y^{(i)}\}$ because we are considering $\max_{x \in \Omega \subset \mathbb{R}^d} f(x)$.

Let denote with $\mathcal{U}$ the set of possible uncertainty quantification measures.

In this paper we consider the following three alternatives:

- GP's predictive standard deviation, namely $\sigma(x)$.
- GP's differential entropy. For a GP it is given by $H(y|X_{1:n}) = \frac{1}{2}\log \det(K) + \frac{d}{2}\log \det(2\pi e)$, where $K \in \mathfrak{R}^{n \times n}$ with entries $K_{ij} = k(x^{(i)}, x^{(j)}), \forall\, x^{(i)}, x^{(j)} \in X_{1:n}$ (Williams & Rasmussen, 2006).
- Distance from previous decisions, inspired from (Bemporad, 2020) and denoted by $z(x)$ :

$$z(x) = \begin{cases} 0 & \text{if } \exists\, x^{(i)} \in \mathbf{X}_{1:n}: \left\| x - x^{(i)} \right\|_2^2 = 0 \\ \frac{2}{\pi} \tan^{-1}\left( \frac{1}{\sum_{j=1}^{n} w_j(x)} \right) & \text{otherwise} \end{cases} \quad (4)$$

with $w_j(x) = \dfrac{e^{-\left\| x - x^{(j)} \right\|_2^2}}{\left\| x - x^{(j)} \right\|_2^2}$.

## 3. Pareto analysis

Given the GP conditioned on the decisions performed so far, it is possible to map the next decision $x^{(n+1)} \in \Omega$ – whichever it is – as a bi-objective choice, with objectives $\zeta(x)$ and $u(x) \in \mathcal{U}$ (both to be maximized).

Pareto rationality is the theoretical framework to analyse multi-objective optimization problems where $q$ objective functions $\gamma_1(x), \ldots, \gamma_q(x)$ where $\gamma_i(x): \to \mathbb{R}$ are to be simultaneously optimized in $\Omega \subseteq \mathbb{R}^d$. We use the notation $\boldsymbol{\gamma}(x) = \left( \gamma_1(x), \ldots, \gamma_q(x) \right)$ to refer to the vector of all objectives evaluated at a location $x$. The goal in multi-objective optimization is to identify the Pareto frontier of $\boldsymbol{\gamma}(x)$.

To do this we need an ordering relation in $\mathbb{R}^q$: $\boldsymbol{\gamma} = (y_1, \ldots, y_q) \preccurlyeq \boldsymbol{\gamma}' = (y_1', \ldots, y_q')$ if and only if $\gamma_i \leq \gamma_i'$ for $i = 1, \ldots, q$. This ordering relation induces an order in $\Omega$: $x \preccurlyeq x'$ if and only if $\boldsymbol{\gamma}(x) \preccurlyeq \boldsymbol{\gamma}(x')$.

We also say that $\gamma'$ dominates $\gamma$ (strongly if $\exists\, i = 1, \ldots, q$ for which $\gamma_i < \gamma_i'$). The optimal non-dominated solutions lay on the so-called Pareto frontier.

The interest in finding locations $x$ having the associated $\boldsymbol{\gamma}(x)$ on the Pareto frontier is clear: they represent the trade-off between conflicting objectives and are the only ones, according to the Pareto rationality, to be considered.

In this paper $q = 2$, with $\gamma_1(x) = \zeta(x)$ and $\gamma_2(x) = u(x) \in \mathcal{U}$. Both the objectives are not expensive to evaluate, therefore the Pareto frontier can be easily approximated by considering a fine grid of locations in $\Omega$ without the need to resort to methods approximating expensive Pareto frontiers within a limited number of evaluations, such as in (Zuluaga et al., 2013).

Thus, we approximate our Pareto frontier by sampling a grid of $m$ points in $\Omega$, denoted by $\widehat{\mathbf{X}}_{1:m} = \left\{ x^{(j)} \right\}_{j=1,\ldots,m}$, and then computing the associated pairs $\Psi_{1:m} = \left\{ \left( \zeta(x^{(j)}), u(x^{(j)}) \right) \right\}_{j=1,\ldots,m}$.

It is important to remark that a large value of $m$ is needed to have a good approximation of the Pareto frontier but this is not an issue because the computational cost is dominated by conditioning the GP on observations (i.e., $\mathcal{O}(n^3)$, with $n \ll m$) instead of making predictions (i.e., inference). The Pareto frontier can be approximated as:

$$\mathcal{P}(\Psi_{1:m}) = \{ \psi \in \Psi_{1:m} : \forall\, \psi' \in \Psi_{1:m}\; \psi \succ \psi' \} \quad (5)$$

where $\psi = \left( \zeta(x), u(x) \right)$ and $\psi' = \left( \zeta(x'), u(x') \right)$, and $\psi \succ \psi' \Leftrightarrow \zeta(x) > \zeta(x') \wedge u(x) > u(x')$.

Figure 1 shows an example of Pareto frontier for $\zeta(x)$ and $u(x) = \sigma(x)$. First five charts, top-left to bottom-right, depict $\Psi_{1:m}$ and the associated $\mathcal{P}(\Psi_{1:m})$ for each kernel in $\mathcal{K} = \{ k_{SE}, k_{EXP}, k_{PE}, k_{M3/2}, k_{M5/2} \}$, separately. The last chart (bottom-right) compares only the five Pareto frontiers, better highlighting the role of the GP kernel. For this example, $f(x)$ is the Branin-Hoo (Jekel & Haftka, 2019) function in $\Omega: [5; 10] \times [0; 15]$, $m = 1976$.

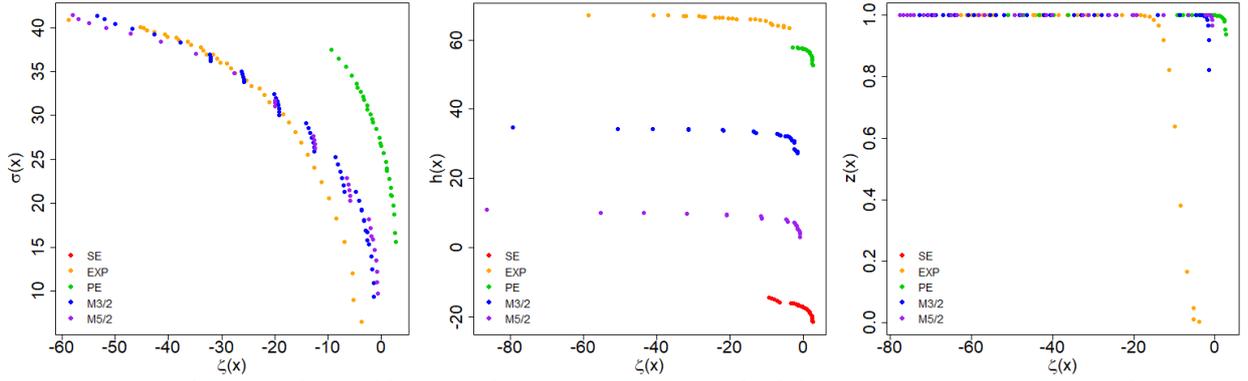

*Figure 1. Pareto frontiers obtained by using the GP's posterior standard deviation as uncertainty quantification (i.e., $u(x) = \sigma(x)$). Five different kernels are used to fit as many GPs, leading to as many Pareto frontiers. Last chart (bottom-right) depicts the five frontiers all together for an easier comparison.*

The only way to analyse how different uncertainty quantification measures can lead to completely different decisions – even if anyway Pareto rational – is to localize, within the search space $\Omega \subset \Re^d$, the locations whose associated objectives lays on the Pareto frontier (namely, the Pareto set). Figure 2 reports just a 2D example, considering ten previous decisions (bold black crosses), five different kernels and three alternative uncertainty quantification measure.

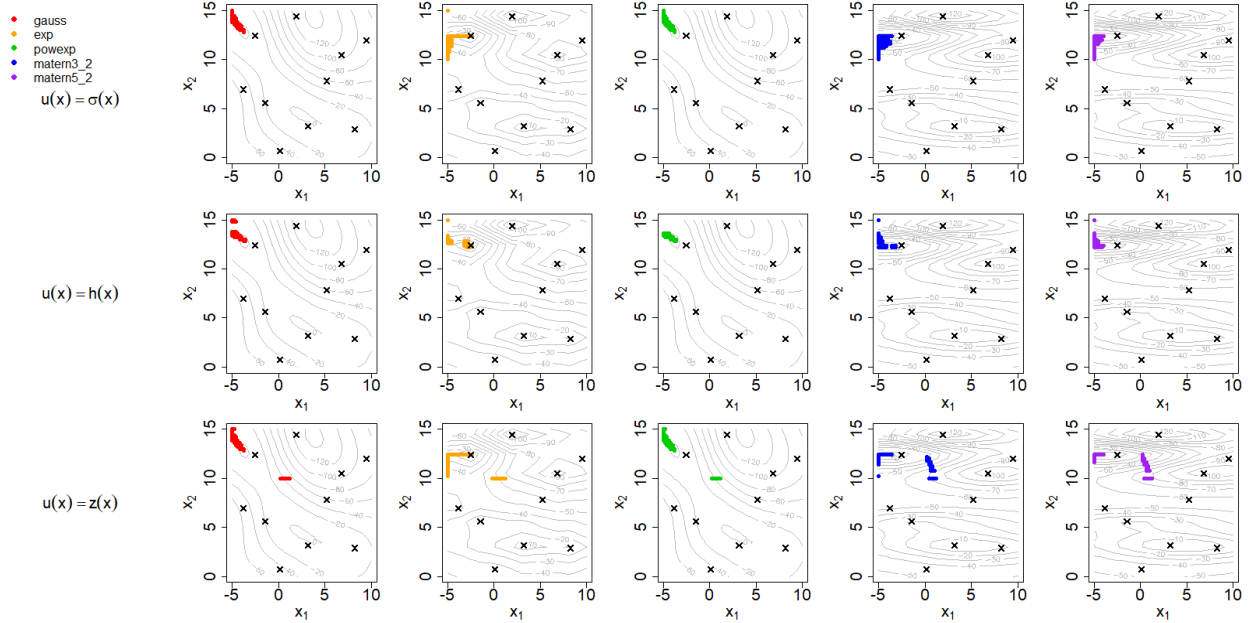

*Figure 2. Next decision depending on: (i) ten previous observations (bold crosses), (ii) uncertainty quantification measures (rows: $\sigma(x)$, $h(x)$ and $z(x)$), and kernels (columns: "gauss" for $k_{SE}$, "exp" for $k_{EXP}$, "powexp" for $k_{PE}$, "matern3_2" for $k_{M3/2}$ and "matern5_2" for $k_{M5/2}$)*

From the figure it is possible to notice that the region of locations associated to Pareto-rational decisions does not change so much depending on kernel, as well as by using $\sigma(x)$ or $h(x)$. The most evident difference arises by using $z(x)$ as uncertainty quantification measure, because it allows to consider as Pareto-rational also decisions in the area around, approximately, the location $(x_1 = 1; x_2 = 11)$. This area is associated to a more explorative behaviour compared to the other – which is also identified by using the other two uncertainty quantification measures – meaning that some explorative choice could be still considered Pareto rational when $u(x) = z(x)$.

This is just an example for explanatory purposes, the hypothesis is investigated and validated in our analysis. Moreover, we have also to consider that humans, (Kahneman, 2011) could take non-Pareto-rational decisions, and it is therefore important to measure how much a decision can be considered *"far from a Pareto-rational one"*. This issue is addressed and formalized in the next section.

**3.1 Distance from the Pareto rationality**

Every next decision, $x^{(n+1)}$, can be analysed according to the distance of its "image" $\left(\zeta(x^{(n+1)}), u(x^{(n+1)})\right)$ from the Pareto frontier, computed as follows:

$$d(\bar{\psi}, \bar{\mathcal{P}}) = \min_{\psi \in \bar{\mathcal{P}}} \left\{ \|\bar{\psi} - \psi\|_2^2 \right\} \tag{6}$$

where $\bar{\psi} = \left(\zeta(x^{(n+1)}), u(x^{(n+1)})\right)$ and $\bar{\mathcal{P}} = \mathcal{P}(\Psi_{1:m}) \cup \{\bar{\psi}\}$.

This distance is computed for every choice among the five kernels and the three uncertainty quantification measures previously presented. The hypothesis is that humans are mostly Pareto rational, and they should therefore make decisions laying on – or close to – the Pareto frontier.
We analyse the distances from all the 15 possible Pareto frontiers (5 kernels × 3 uncertainty quantification measures) and how they change along the optimization process. Figure 3 shows an example taken from our experimental results and anticipated here just for explanatory purposes. The 10 charts refer to as many black box optimization problems solved by a single human subject. At each iteration, and for each uncertainty quantification measure, the minimum distance between the associated Pareto frontier and the human decision is reported, irrespectively to the kernel. It is important to remark that distances cannot be compared in absolute terms, because the three different uncertainty measures can vary in very different ranges. However, it is possible to observe that distances result correlated in some cases and uncorrelated in others.
Finally, from the charts it is possible to notice that: *(a)* in some cases Pareto rationality is independent on the uncertainty quantification, such as for the problems: *bukin6*, *goldpr*, *rastr*, *stybtang*, but not in general; *(b)* a higher number of decisions are considered Pareto rational if $u(x) = z(x)$; *(c)* in some cases it is possible to observe a shift from Pareto-rationality to *not*-Pareto rationality (e.g., this is evident in for *beale, goldpr* and *rastr*).

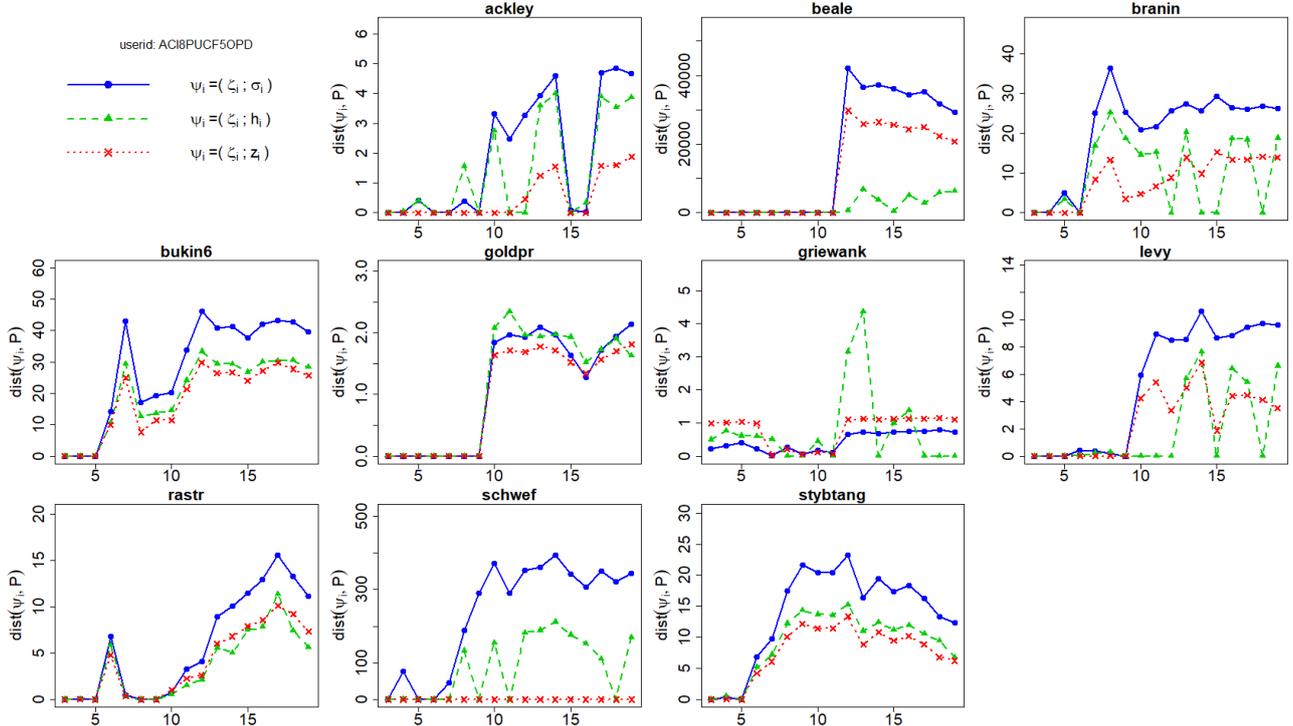

*Figure 3. Distance, at each iteration, of the next decision from three different Pareto frontiers, one for each uncertainty quantification measure $\sigma(x)$, $h(x)$ and $z(x)$. All the fifteen charts are related to as many black box optimization tasks performed by a single human subject.*

## 4. The Wasserstein distance – Basic notions and numerical approximation

Measuring the distance between distributions can be accomplished by many alternative models. A general class of distances, known as f-divergences, is based on the expected value of a convex function of the ratio of two distributions. If $P$ and $Q$ are two probability distributions over $\mathbb{R}^d$ and $f$ is a convex function such that $f(0) = 1$ the f-divergence is given by:

$$D_f(P,Q) = \mathbb{E}_Q f\left(\frac{P}{Q}\right) \tag{7}$$

According to the choice of $f$ the above formula yields specific distances including Kullback-Leibler (and its symmetrized version Jensen-Shannon), Hellinger, total variation and $\chi$-square divergence.
In this paper we focus on the Wasserstein distance whose basic notions are given in Sect. 4.1 while Sect. 4.2 is devoted to the computation of the barycenter between distributions and the extension of k-means clustering to the Wasserstein space. It is important to remark that the presentation is quite basic omitting any mathematical characterization of WST for which the reader is referred to (Villani, 2009) and (Peyre & Cuturi, 2019). The WST metric is based on the solution of an optimal transport problems. WST enables to synthetizes the comparison between two multi-dimensional distributions through a single metric using all information in the distributions. Moreover, WST is generally well defined and provide an interpretable distance metric between distributions.
The Wasserstein distance can be traced back to the works of Gaspard Monge (1781) and Lev Kantorovich (1942). Recently, also under the name of Earth Mover Distance (EMD) it has been gaining increasing importance in several fields like Imaging (Bonneel et al., 2016), Natural Language Processing (NLP) (Huang et al., 2016) and the generation of adversarial networks (Arjovsky et al., 2017). Important references are (Peyré & Cuturi, 2019; Villani, 2009), which also gives an up-to-date survey of numerical methods

**4.1 Basic notions**

The WST distance between continuous probability distributions is:

$$W_p\big(P^{(1)}, P^{(2)}\big) = \left(\inf_{\gamma \in \Gamma(P^{(1)}, P^{(2)})} \int_{X \times X} d\big(x^{(1)}, x^{(2)}\big)^p d\gamma\big(x^{(1)}, x^{(2)}\big)\right)^{\frac{1}{p}} \quad (8)$$

where $d(x^{(1)}, x^{(2)})$ is also called *ground distance* (usually it is the Euclidean norm), $\Gamma(P^{(1)}, P^{(2)})$ denotes the set of all joint distributions $\gamma(x^{(1)}, x^{(2)})$ whose marginals are respectively $P^{(1)}$ and $P^{(2)}$, and $p > 1$ is an index. The Wasserstein distance is also called the Earth Mover Distance (EMD). Intuitively, Figure 4 indicates how much mass must be transported from $x^{(1)}$ to $x^{(2)}$ to transform the distributions $P^{(1)}$ into the distribution $P^{(2)}$. The EMD is the minimum energy cost of moving and transforming a pile of sand in the shape of $P^{(1)}$ to the shape of $P^{(2)}$. The cost is quantified by the amount of sand moved times the moving distance $d(x^{(1)}, x^{(2)})$. The EMD then is the cost of the optimal transport plan.

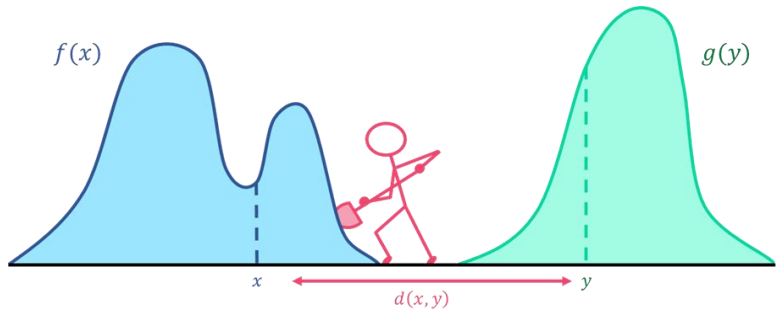

*Figure 4. Behind the Earth Mover Distance*

Figure 5 represents the traditional point-wise distance between distributions.

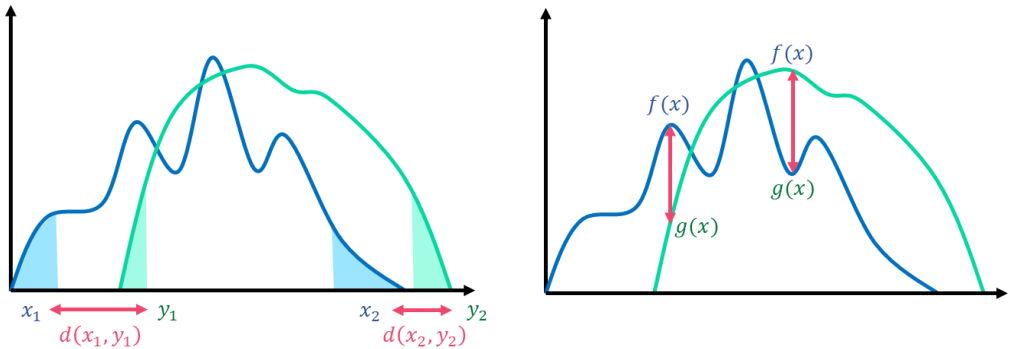

*Figure 5. Point-wise distance between two distributions*

There are some specific cases, very relevant in applications, where WST can be written in an explicit form. Let $\hat{P}^{(1)}$ and $\hat{P}^{(2)}$ be the cumulative distribution for one-dimensional distributions $P^{(1)}$ and $P^{(2)}$ on the real line and $(\hat{P}^{(1)})^{-1}$ and $(\hat{P}^{(2)})^{-1}$ be their quantile functions.

$$W_p(P^{(1)}, P^{(2)}) = \left( \int_0^1 \left| (\hat{P}^{(1)})^{-1}(x^{(1)}) - (\hat{P}^{(2)})^{-1}(x^{(2)}) \right|^p dx \right)^{\frac{1}{p}} \tag{9}$$

Let's now consider the case of a discrete distribution $P$ specified by a set of support points $x_i$ with $i = 1, \ldots, m$ and their associated probabilities $w_i$ such that $\sum_{i=1}^m w_i = 1$ with $w_i \geq 0$ and $x_i \in M$ for $i = 1, \ldots, m$. Usually, $M = \mathbb{R}^d$ is the $d$-dimensional Euclidean space with the $l_p$ norm and $x_i$ are called the support vectors. $M$ can also be a symbolic set provided with a symbol-to-symbol similarity. $P$ can also be written using the notation:

$$P(x) = \sum_{i=1}^m w_i \delta(x - x_i) \tag{10}$$

where $\delta(\cdot)$ is the Kronecker delta.

The WST distance between two distributions $P^{(1)} = \{w_i^{(1)}, x_i^{(1)}\}$ with $i = 1, \ldots, m_1$ and $P^{(2)} = \{w_i^{(2)}, x_i^{(2)}\}$ with $i = 1, \ldots, m_2$ is obtained by solving the following linear program:

$$W(P^{(1)}, P^{(2)}) = \min_{\gamma_{ij} \in \mathbb{R}^+} \sum_{i \in I_1, j \in I_2} \gamma_{ij} \, d\left( x_i^{(1)}, x_j^{(2)} \right) \tag{11}$$

The cost of transport between $x_i^{(1)}$ and $x_j^{(2)}$, $d\left( x_i^{(1)}, x_j^{(2)} \right)$, is defined by the $p$-th power of the norm $\left\| x_i^{(1)}, x_j^{(2)} \right\|$ (usually the Euclidean distance).
We define two index sets $I_1 = \{1, \ldots, m_1\}$ and $I_2$ likewise, such that

$$\sum_{i \in I_1} \gamma_{ij} = w_j^{(2)}, \forall j \in I_2 \tag{12}$$

$$\sum_{j \in I_2} \gamma_{ij} = w_i^{(\ )}, \forall i \in I_1 \tag{13}$$

Equations 12 and 13 represent the in-flow and out-flow constraint, respectively. The terms $\gamma_{ij}$ are called matching weights between support points $x_i^{(1)}$ and $x_j^{(2)}$ or the optimal coupling for $P^{(1)}$ and $P^{(2)}$.
The discrete version of the WST distance is usually called Earth Mover Distance (EMD). For instance, when measuring the distance between grey scale images, the histogram weights are given by the pixel values and the coordinates by the pixel positions. Another way to look at the computation of the EMD is as a network flow problem. In the specific case of histograms, the entries $\gamma_{ij}$ denote how much of the bin $i$ has to be moved to bin $j$.
The basic computation of OT between 2 discrete distributions involves solving a network flow problem whose computation scales typically cubic in the sizes of the measure. There are 2 lines of work to reduce the time complexity of OT (i) simple ground costs can lead to simpler computations. In the general case it is shows to be equivalent to a min-flow algorithm of quadratic computational complexity and, in specific cases, to linear.

The computation of EMD turns out to be the solution of a minimum cost flow problem on a bi-partite graph where the bins of $P^{(1)}$ are the source nodes and the bins of $P^{(2)}$ are the sinks while the edges between sources and sinks are the transportation costs.

In the case of one-dimensional histograms, the computation of WST reduces to the comparison of two 1-dimensional histograms which can be performed by a simple sorting and the application of the following equation (14).

$$W_p\left(P^{(1)}, P^{(2)}\right) = \left(\frac{1}{n}\sum_{i}^{n}\left|x_i^{(1)*} - x_i^{(2)*}\right|^p\right)^{\frac{1}{p}} \tag{14}$$

where $x_i^{(1)*}$ and $x_i^{(2)*}$ are the sorted samples.

**4.2 Barycenter and clustering**

Consider a set of $N$ discrete distributions, $\mathbf{P} = \{P^{(1)}, \dots, P^{(N)}\}$, with $P^{(k)} = \left\{\left(w_i^{(k)}, x_i^{(k)}\right) : i = 1, \dots, m_k\right\}$ and $k = 1, \dots, N$, then, the associated barycenter, denoted with $\bar{P} = \{(\bar{w}_1, x_1), \dots, (\bar{w}_m, x_m)\}$, is computed as follows:

$$\bar{P} = \underset{P}{\text{argmin}} \frac{1}{N}\sum_{k=1}^{N}\lambda_k W\left(P, P^{(k)}\right) \tag{15}$$

where the values $\lambda_k$ are used to weight the different contributions of each distribution in the computation. Without loss of generality, they can be set to $\lambda_k = 1/N \ \forall \ k = 1, \dots, N$.

Among the advantages of WST, we remark:
- It is sensitive to the underlying geometry. Consider 3 distributions $P^{(1)} = \delta_0, P^{(2)} = \delta_\varepsilon$ and $P^{(3)} = \delta_{100}$. $W(P^{(1)}, P^{(2)}) \approx 0$, $W(P^{(1)}, P^{(3)}) \approx W(P^{(2)}, P^{(3)}) \approx 100$. The distances Total variation, Hellinger and Kullback-Leibler take the value 1, thus they fail to capture our intuition that $P^{(1)}$ and $P^{(2)}$ are close to each other while they are far away from $P^{(3)}$.
- It is *shape preserving*. Denote $P^{(1)}, \dots, P^{(N)}$ and assume that each $P^{(j)}$ can be written as a location shift of any other $P^{(i)}$, with $i \neq j$. Suppose that each $P^{(j)}$ is defined as $P^{(j)} = \mathcal{N}(\mu_j, \Sigma)$, then the barycenter has the closed form:

$$\bar{P} = \mathcal{N}\left(\frac{1}{N}\sum_{j=1}^{N}\mu_j, \Sigma\right)$$

in contrast to the (Euclidean) average of the $\frac{1}{N}\sum_{j=1}^{N}P^{(j)}$.

Therefore, the concept of barycenter enables clustering among distributions, in a space whose metric is WST. More simply, the barycenter in a space of distributions is the analogue of the centroid in a Euclidean space. The most common and well-known algorithm for clustering data in the Euclidean space is k-means. Since it is an iterative distance-based (aka representative based) algorithm, it was easy to propose variants of k-means by simply changing the distance adopted to create clusters, such as the Manhattan distance (leading to k-medoids) or any kernel allowing for non-spherical clusters (i.e., kernel k-means). The crucial point is that only distance is changed, while the overall iterative two-step algorithm is maintained. This is also valid in the case of the WST k-means, where the Euclidean distance is replaced by WST and centroids are replaced by barycenters:

- **Step 1 – Assign.** Given the current $k$ barycenters at iteration $t$, namely $\bar{P}_t^{(1)}, \ldots, \bar{P}_t^{(k)}$, clusters $C_t^{(1)}, \ldots, C_t^{(k)}$ are identified by assigning each one of the distributions $P^{(1)}, \ldots, P^{(N)}$ to the closest barycenter:

$$C_t^{(i)} = \left\{ P^{(j)} \in \mathbf{P} : \bar{P}_t^{(i)} = \operatorname*{argmin}_{Q=\{\bar{P}_t^{(1)}, \ldots, \bar{P}_t^{(k)}\}} W(Q, P^{(j)}) \right\}, \forall\, i = 1, \ldots, k$$

- **Step 2 – Optimize.** Given the updated composition of the clusters, update the barycenters:

$$\bar{P}_{t+1}^{(i)} = \operatorname*{argmin}_{Q} \frac{1}{\left|C_t^{(i)}\right|} \sum_{P \in C_t^{(i)}} W(Q, P)$$

that comes directly from equation (15).

As in k-means, a key point of WST k-means is the initialization of the barycenters. In the case that all the distributions in **P** are defined on the same support, then they can be randomly initialized, otherwise, a possibility is to start from $k$ distributions randomly chosen among those in **P**. Finally, termination of the iterative procedure occurs when the result of the assignment step does not change any longer or a prefixed maximum number of iterations is achieved.

The major computational issue is the polynomial complexity of the linear programming solvers commonly used to compute WST. Starting from the consideration that the variables in $w$ are more important that the matching weights, approximate solvers have been proposed, specifically Sinkhorn solvers, which will be detailed later. Here it is just important to remark that they allow to manage the trade-off between accuracy and computational cost through a regularization hyperparameter. Another approach is taken in (Ye et al., 2017) based on ADMM. Entropic regularization enables scalable computations, but large values of the regularization parameter could induce an undesirable smoothing effect while low values not only reduce the scalability but might induce several numeric instabilities.

## 5. Experimental setup

### 5.1. Data structure

All the results are stored into a data table with columns: *user_id*, *problem_id*, *n+1*, *uncertainty_measure*, *min_dist_from_Pareto_frontier*. Finally, numbers of Pareto rational decisions are separately computed for each uncertainty quantification measure and aggregated by *(a)* players and *(b)* test problems.
All the data have been organized into a dataset having the following features, as computed from all the games performed by the players:
- **tf** – test function
- **user** – the user id
- **iter** – iteration (the next decision to take)
- **uq** – uncertainty quantification measure
- **dst** – distance from the Pareto front (having **uq** as uncertainty quantification measure)
- **cum.reward** – the score accumulated so far
- **class** – the type of decision taken by the player (i.e., Pareto or not-Pareto).

All the data can be found at the github repository https://github.com/acandelieri/humans_strategies_analysis

## 5.2 Data analysis

The main goal of data analysis is twofold: to analyze the decisions performed by every volunteer, to characterize their aggregation and to analyse how ACR can explain non Paretian behaviour and how it can be explained in terms of uncertainty quantification and Pareto rationality.
The data collection is organized in 2 consecutive steps:

- *Step 1:* computing the number of Pareto rational decisions, depending on the uncertainty quantification measure, and comparing them. A decision is considered Pareto rational if the distance from the Pareto frontier is less than 0.5. This analysis step is summarized as follows:
    For each player and each test problem do:
    - Initialize $n = 3$ (i.e., the first three decisions of each user and for each test problem cannot be analysed, because fitting a GP over a 2-dimensional search space requires at least three observations).
    - Condition a GP for each one the kernels in $\mathcal{K}$ to the previous $n$ decisions and observations performed by that player for that test problem.
    - For each $u(x) \in \mathcal{U}$ do:
        - Use each one of the conditioned GPs to approximate the associated Pareto frontier by sampling a grid of $m = 30 \times 30 = 900$ locations in $\Omega$ and then computing the associated locations in the $\zeta$-$u$ space.
        - Map $x^{(n+1)}$ into five associated locations $\psi = \left(\zeta(x^{(n+1)}), u(x^{(n+1)})\right)$ – one for each kernel – and compute its minimum distance from the five Pareto frontiers and store it.
        - $n \leftarrow n + 1$ and go to Step 3
    - End For each $u(x) \in \mathcal{U}$
    End For each player and each test problem

- *Step 2:* depending on results from the two previous steps, the uncertainty quantification measure which allows to more frequently classify the humans' choices as Pareto rational is selected. Then, the relationship between the fact that the decision is Pareto rational, and the *reward* collected so far by the user is investigated, with the aim to identify a possible motivation for *over-explorative* decisions (i.e., not-Pareto rational decisions). In our analysis *reward* is represented by the score immediately observed by the player implied by his/her own decision. *Cumulative reward* is therefore the sum of scores collected up to a given decision. Finally, the *Average Cumulated Reward* (ACR), up to a given decision, is computed as the arithmetic mean of the cumulated reward up to that decision:

$$ACR^{(n+1)} = \frac{1}{n}\sum_{i=1}^{n} y^{(i)}$$

we use $ACR^{(n+1)}$ to denote the average reward collected up to $n$ just to be coherent with indexing, since this value is analysed in relation with the Pareto distance of decision $x^{(n+1)}$. The idea is that ACR could quantify the amount of "gratification" (high values of ACR) or "stress" (low values of ACR) experienced by the player in solving the test problem. The hypothesis is that lower values of ACR could be associated to not-Pareto rational decisions, induced by a sense of stress for the incapability to (further) improve the score.

# 6. Experimental results and their analysis

## 6.1 Wasserstein analysis of the test functions

The histogram we use is based on the notion of decile. A decile rank arranges the data in order from lowest to highest and is done on a scale of one to ten where each successive decile corresponds to an increase of 10 percentage points. The basic histogram therefore has 10 bins corresponding to deciles in the distributions whose weights represent the number of players whose decisions were Paretian with that decile.

A stacked histogram is provided for each test problem, comparing the distributions obtained for each one of the three uncertainty quantification measures. The full set of charts is reported in the supplementary material. A first insight can be performed by visual inspection. Clearly *stybang* and *bukin6* are difficult in that they generate fewer Paretian decisions than *schwefel*. If we look at the impact of the uncertainty measure we note $u(x) = z(x)$ is associated to an increase in terms of number of Pareto rational decisions, as shown, for instance, for the test problems *schwef* and *ackley*: for 3 players, 10% of their decisions were Pareto rational when $u(x) = \sigma(x)$ – and no one for the other two uncertainty quantification measures. The number of players substantially increases with the percentage of Pareto rational decisions only for $u(x) = z(x)$. This is taken to the extreme in the case of the *schwef* test function.

The ideal distribution, that is a fully Pareto compliant distribution, is a useful target which enables an intuitive yet quantitative evaluation of the "Pareto value" of each histogram as represented in Figure 6 through the WST distance between each histogram and the ideal one. A further WST based analysis can be conducted according to the distance between a function or a subject histogram and its Paretian ideal. This is reported in table 1.

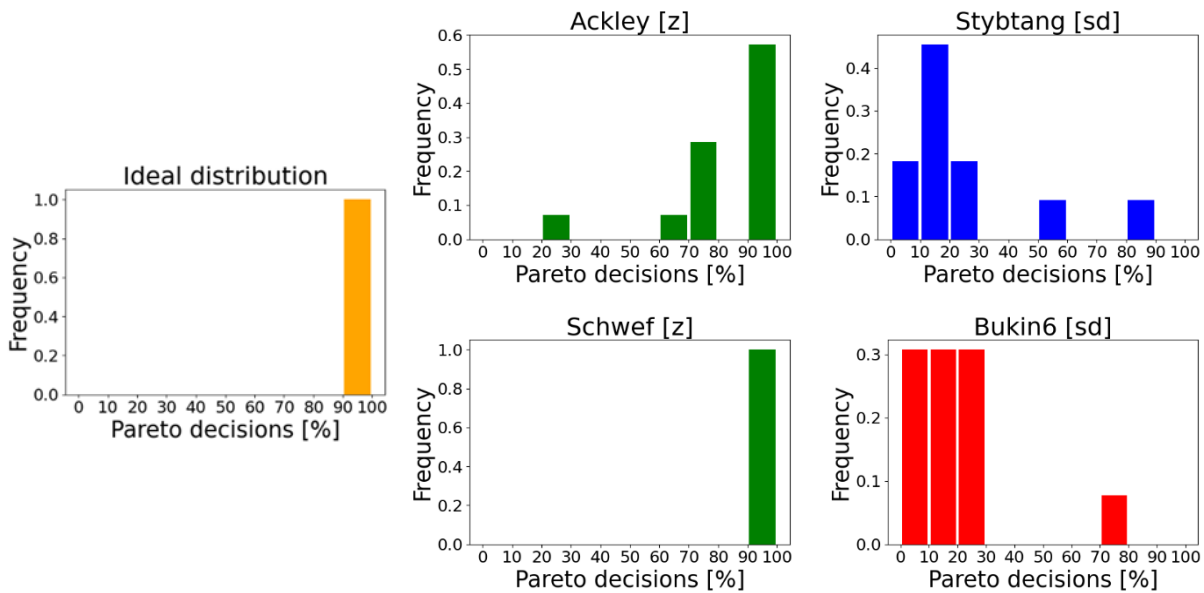

*Figure 6. Number of players with respect to percentage of decisions classified as Pareto rational, separately for the three uncertainty quantification measures. One chart for each test problem.*

*Table 1. Wasserstein distances from the ideal distribution.*

| Test function | H | SD | Z |
|---|---|---|---|
| ackley | 3.57 | 4.79 | **1.29** |
| beale | **6.77** | 7.23 | 7.15 |
| branin | **6.43** | 7.73 | 6.92 |
| bukin6 | **7.23** | 7.54 | 7.25 |
| goldpr | 6.21 | **4.36** | 5.50 |
| griewank | 4.00 | **0.57** | 3.86 |
| levy | 4.79 | 6.43 | **2.86** |
| rastr | **3.64** | 4.21 | 3.79 |
| schwef | 4.71 | 5.23 | **0.00** |
| stybtang | **6.92** | 7.00 | 7.08 |
| **Barycenter** | 5.33 | 5.52 | **4.58** |

The data can be evaluated globally by using the barycenter calculated according to formulas (15) and shown in Figure 7 (the distance from the barycenter to the ideal histogram is in the last row of Table 1).

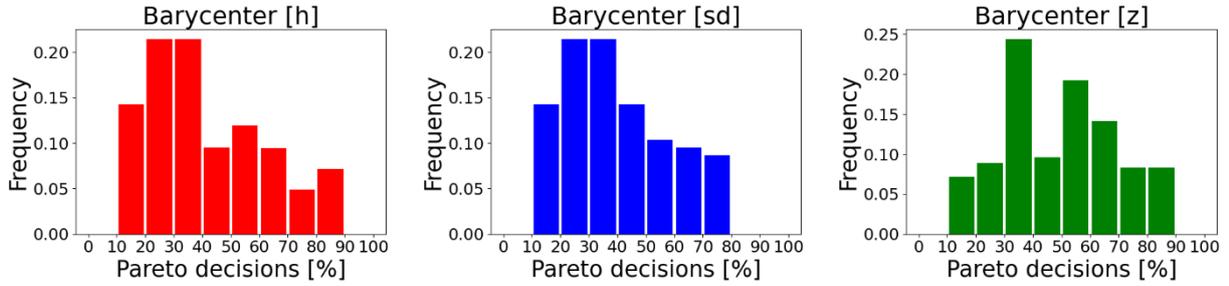

*Figure 7. Wasserstein barycenters of the histograms related to the test functions.*

Also clustering can be performed in the WST space. Specifically, we have used WST k-means. Since our main objective is to partition the behavioural patterns into Pareto and not-Pareto decisions, k=2 is a reasonable choice. The results depend on the uncertainty measure and are:

In the case $u(x) = h(x)$ clusters are:
$$C_1^h = (griewank, levy, rastr, schwef, ackley),$$
$$C_2^h = (stybtang, goldpr, beale, bukin6, branin)$$

In the case $u(x) = \sigma(x)$ clusters are:
$$C_1^{sd} = (griewank, rastr, schwef, ackley, goldpr)$$
$$C_2^{sd} = (stybtang, levy, beale, bukin6, branin)$$

In the case $u(x) = z(x)$ clusters are:
$$C_1^z = (griewank, rastr, schwef, ackley, levy)$$
$$C_2^z = (stybtang, goldpr, beale, bukin6, branin)$$

The functions which were visually singled out as "hard" and "easy" are correctly assigned to two different clusters under the uncertainty measure $h(x)$ and $z(x)$. We have computed also $k = 3$ and $k = 4$, the interpretation is less natural, and the metric value accordingly worsens.

The values in Table 1 of the distance from the ideal are summarized if Figure 8. Two relevant results have to be remarked: (i) the two identified clusters always capture also the difference in terms of distance from ideal independently of the uncertainty quantification measures, and (ii) the distinction between the two clusters is maximized for the uncertainty quantification measure $z(x)$.

A way of visualizing the clustering results is through the box-plot representation, which shows that cluster 1, associated to Pareto decisions, has a significantly smaller WST distance than cluster 2, even more significant in the case that z is used as uncertainty quantification measure.

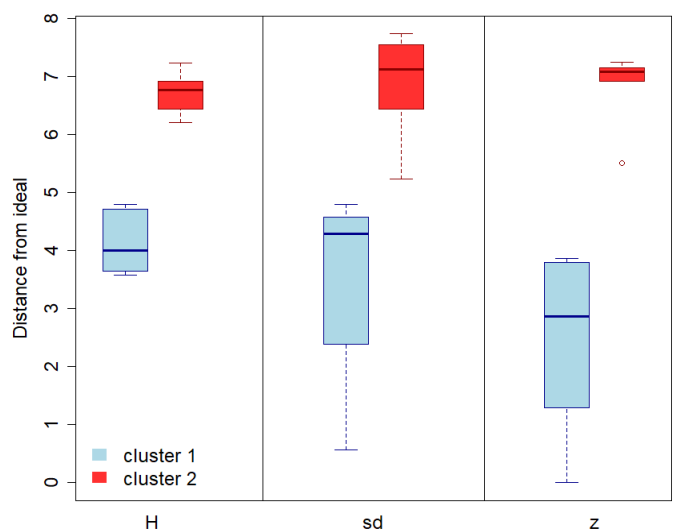

**Figure 8.** *WST distance from the ideal, represented as boxplots, for the two clusters and, separately, for each uncertainty quantification measure.*

## 6.2 Wasserstein analysis of the users

Figure 9 shows histograms like those reported in the previous section but referred to each subject: the weights are the distribution of the number of test problems with respect to the percentage of Pareto rational decisions. Again, we can notice, by visual inspection, two relatively Paretian players (U01 and U13) and two not-Paretian (U05 and U14); moreover, the highest percentages of Pareto rational decisions are associated to $z(x)$.

The ideal distribution, that is a fully Pareto compliant distribution, is useful target which enables an intuitive yet quantitative evaluation of the "Pareto value" of each histogram as represented in Table 2 through the WST distance between each histogram and the ideal one.

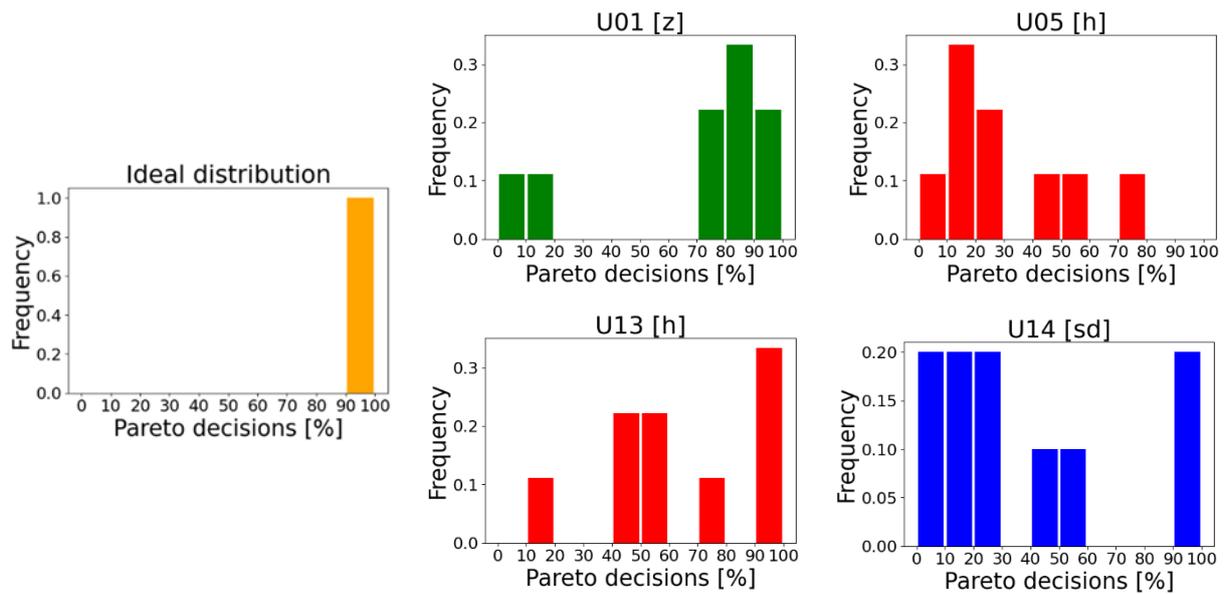

*Figure 9.* Number of test problems with respect to percentage of decisions classified as Pareto rational, separately for the three uncertainty quantification measures. One chart for each player.

*Table 2.* Wasserstein distances from the ideal distribution.

| User | H | SD | Z |
| --- | --- | --- | --- |
| U01 | 3.89 | 4.11 | **2.67** |
| U02 | 5.50 | 5.22 | **4.89** |
| U03 | 6.10 | **5.00** | 5.10 |
| U04 | 5.80 | 5.90 | **4.56** |
| U05 | 6.44 | 5.89 | **5.00** |
| U06 | 5.30 | 5.20 | **3.00** |
| U07 | 5.30 | 5.78 | **5.00** |
| U08 | 5.50 | 5.40 | **3.80** |
| U09 | 5.33 | **4.14** | 4.50 |
| U10 | 5.80 | 5.60 | **4.78** |
| U11 | 5.60 | 6.10 | **4.30** |
| U12 | **5.80** | 6.80 | **5.80** |
| U13 | **3.11** | 4.20 | 3.50 |
| U14 | 5.60 | 5.70 | **4.78** |
| **Barycenter** | 5.47 | 5.34 | **4.40** |

This effect can be evaluated globally using the barycenter (Figure 10) computed according to the formula and the distance from the barycenter to the ideal situation (last row of Table 2).

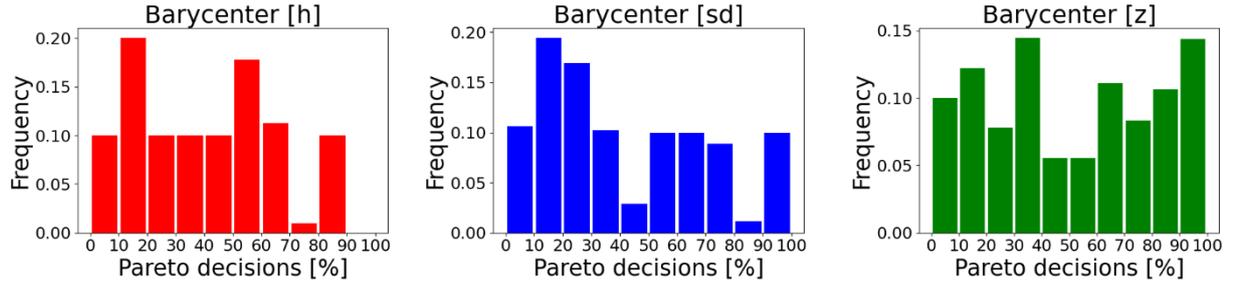

*Figure 10. Wasserstein barycenters of the histograms related to the users.*

Since our main objective is to partition the behavioural patterns into Pareto and not-Pareto decisions, $k = 2$ seems a reasonable choice for WST k-means clustering. As a result, the two clusters can capture most of the difference in terms of distance from ideal, independently of the adopted uncertainty quantification measure, but now the difference is more relevant in the case of the entropy-based uncertainty quantification measure, instead of the inverse distance based one.

In the case $u(x) = h(x)$ clusters are:
$$C_1^h = (U01, U13)$$
$$C_2^h = (U02, U03, U04, U05, U06, U07, U08, U09, U10, U11, U12, U14)$$

In the case $u(x) = \sigma(x)$ clusters are:
$$C_1^{sd} = (U01, U03, U09, U13),$$
$$C_2^{sd} = (U02, U04, U05, U06, U07, U08, U10, U11, U12, U14);$$

In the case $u(x) = z(x)$ clusters are:
$$C_1^z = (U01, U06, U08, U09, U11, U13),$$
$$C_2^z = (U02, U03, U04, U05, U07, U10, U12, U14)$$

The subjects which were visually singled out as Paretian and non-Paretian are separately assigned to two clusters. We have computed also $k = 3$ and $k = 4$, the interpretation is less natural, and the metric value accordingly worsens. A way of visualizing the clustering results is through the box-plot representation which shows that cluster 1, associated to Pareto decisions, has a significantly smaller WST distance than cluster 2.

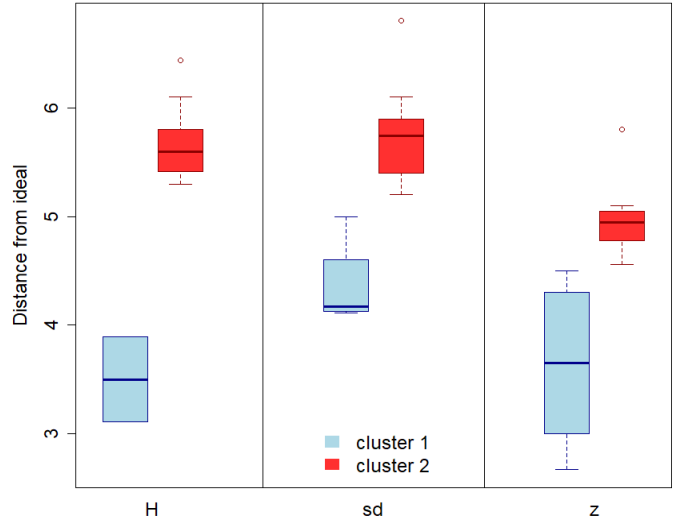

*Figure 11.* WST distance from the ideal, represented as boxplots, for the two clusters and, separately, for each uncertainty quantification measure.

### 6.3 Explaining Pareto and not-Pareto decisions via a Decision Tree classifier

Now we take up specifically the issue of non-Paretian behaviour. The previous results show that a substantial number of decisions show a significant deviation from Pareto. One explanation is that the perceived impossibility to further improve with respect to the best score obtained might generate a level of stress and lead players to the so-called "exasperated exploration". This non-Pareto behaviour should be captured by the Average Cumulated Reward (ACR) in (Candelieri et al., 2021).

In this section we investigate, more deeply, the possible reasons underlying the "exasperated exploration" behaviour – which can be basically assimilated to Pure Random Search. As we want that these reasons could be interpretable, we have decided to use a Decision Tree classifier (i.e., J48 from the suite Weka) (Hall et al., 2009). The dataset consists of all the choices performed by players along their games (excluded the first three "shots" for each player and for each game, because they are obviously explorative). Thus, each row (aka instance) of the dataset is a decision. The columns (aka features) of the dataset are: (1) the black-box function underlying the game -namely *tf*, (2) the user identifier – namely *user*, (3) the iteration at which the decision has been taken-namely *iter*, (4) the cumulated reward - namely *cum.reward*, and finally (5) the type of decision, Pareto/not-Pareto, that is the "class label" to be predicted by the Decision Tree.

Figure 12 reports the learned DT. The relevant considerations from this model are:

1. The feature that has mostly characterized the players' decisions in our study is the underlying black-box function (*tf*). This means that some problems are "easy" and encourage/facilitate Pareto decisions (i.e., *ackley*, *griewank*, *levy* and *schwef*), independently on the player, iterations, and cumulated reward, while others are more complicated and lead to "exasperated exploration" as main behaviour to approach them (i.e., *beale*, *branin*, *bukin*, *goldpr*, *stybtang*). This conclusion matches the results of previous section.
2. Distinguishing between Pareto and not-Pareto decisions in the case of the *rastr* problem requires to consider other two features, that are *user* and *iter*. Indeed, it is reasonable to think that decisions can depend on both the personal risk-attitude (*user*) and perceived sense of "stress". More precisely, decisions performed by 5 out of the 14 players are classified as Pareto by the learned DT, while for the remaining 9 players the feature *iter* becomes crucial: as soon as iterations exceed a certain value, decisions are classified as not-Pareto. The value to exceed is different among the 9 players, depending either on a different perception of the "stress" or the initial set of decisions performed by the player which has implied a completely different score collection along his/her search process. In any case, the rule related to *iter* remarks that not-Pareto decisions increase with approaching the end of the game.

Finally, just to avoid that order of the features could affect the learned model, we have inverted the order of the features (except for the class) and performed again the analysis with J48. In this case the best value for the J48's confidence factor resulted 0.02.
As a result, we have obtained a DT having training accuracy=78.36% and validation accuracy=78.23%. Accuracies are in line with the previous ones, just a slight improvement in validation accuracy is observed.

Figure 13 reports the learned DT. The relevant considerations from this model are:
1. As for the previous DT, (the same) 4 test problems are considered as "easy" and implying Pareto decisions, while the number of complicated test problems decreases from 5 to 4 (*goldpr* is no more in this set).
2. For the test problem *goldpr*, the feature *user*, in case combined with *cum.reward* or *iter*, allows to distinguish between Pareto and not-Pareto decisions. Analogously to *iter*, also *cum.reward* is an indicator of possible stress: the lower the *cum.reward* the higher the stress.
3. More interestingly, *cum.reward* enables simpler classification rule with respect to the previous DT. More precisely, the level of "stress" – quantified by *cum.reward* – is more relevant than *user* (which is not considered any longer). This result is in line with those related to ACR and reported in (Humans-II).

We report a simplified and more "compact" representation of the DT. More specifically we have collapsed branches related to individual users having the same predicted decision into just one single branch, still indicating the user-dependence of thresholds on that collapsed branch.

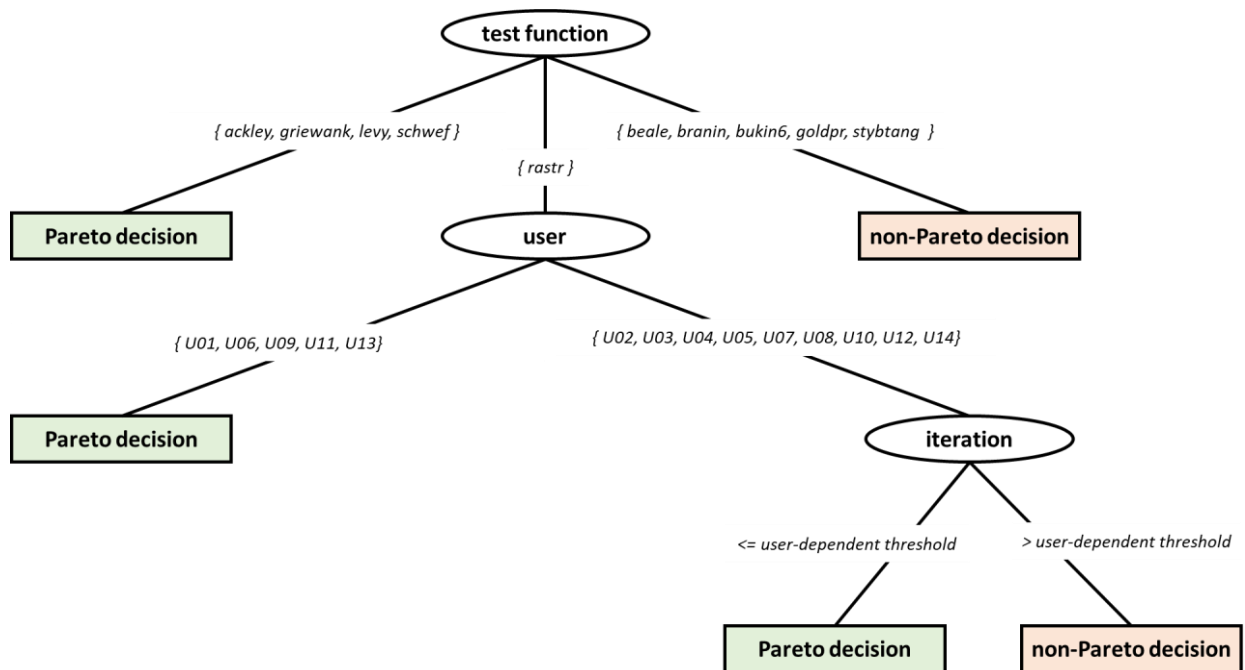

*Figure 12. Decision tree to classify between Pareto and not-Pareto decisions.*

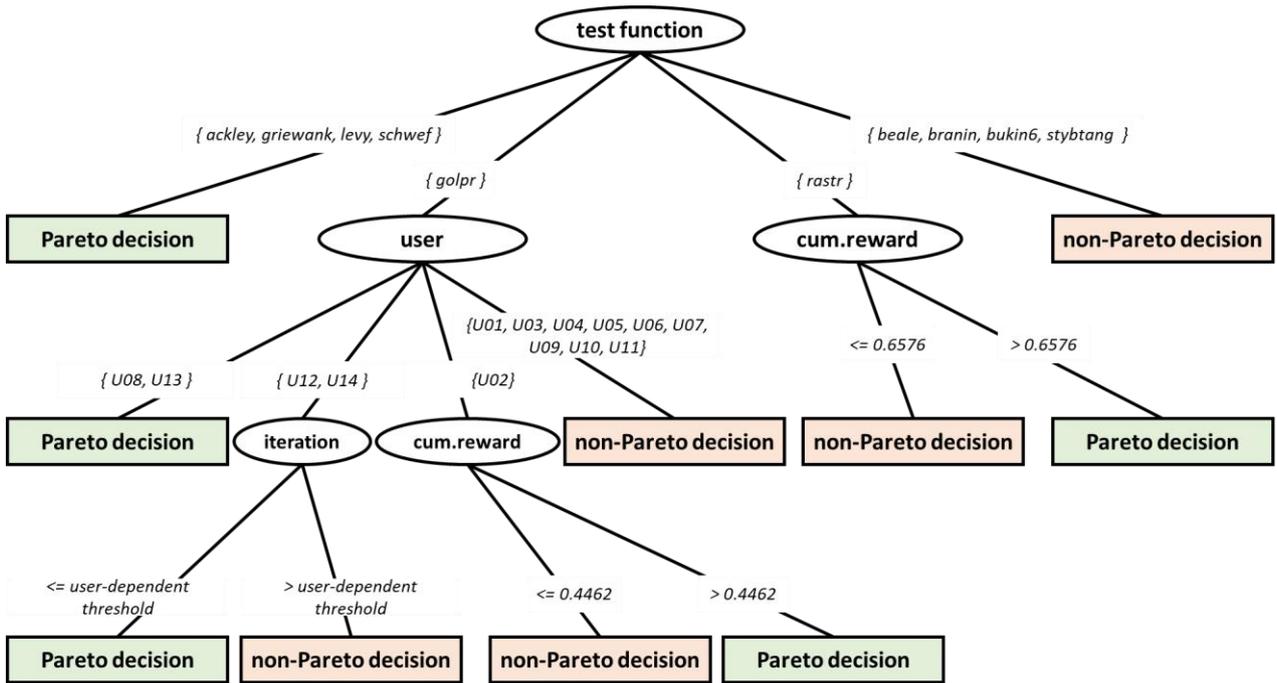
*Figure 13. Decision tree to classify between Pareto and not-Pareto decisions (with features resorted).*

## 7. Conclusions and perspectives: learning with distances

The Pareto analysis of human search data in the Wasserstein space has been shown to be a useful representation of how different uncertainty quantifications relate to human behaviour when solving optimization task.
The distributional analysis in the WST space has been shown to be a useful tool to yield new insights to characterize the human behaviour both at the individual level (single user/single task) and an aggregate level computing barycenters and performing clustering in the Wasserstein space.

According to a Decision Tree based analysis of the human search data, cumulative reward collected over the search process and number of decisions performed so far resulted be useful indicators for explaining the onset of "exasperated exploration", associated to not-Pareto decisions. This result was consistent with those obtained from the WST-based analysis.

A main motivation in writing this paper stems from the observation that *human learners* are fast and effective at managing the exploration/exploitation dilemma, processing the upcoming knowledge, and adapting to unfamiliar environments: the characterization of this behaviour, besides being a key question for cognitive and behavioural sciences, is also an important opportunity for Machine Learning.
There are two main questions underlying the results of this paper: (i) might we learn from gamification experiments something about how humans take to uncertainty searching for a goal? (ii) Might this offer new insights into the design of algorithms of machine learning and optimization? These points can be better specified by the following research questions:

- Do humans always make "rational" choices (i.e., Pareto optimal decisions in the space $\psi = (\zeta(x), u(x))$ or, in some cases, they "exasperate" exploration? The statistical analysis of the results does not disprove that most of the decisions is Pareto rational.
- A further step should be a probabilistic characterization of the sequence of decisions and a close analysis of the dynamics of how people change their behaviour.

- Do different uncertainty quantification measures lead to different classifications of humans' decisions? And which uncertainty quantification measure make humans "more rational"?
- Another result is that humans do not always make "rationale" choices (i.e., Pareto optimal decisions in the space of expected improvement and uncertainty) and in some cases, they "exasperate" exploration. The analysis of computational results allows to formulate at least a tentative answer to why or rather in which conditions we observe deviations from "rationality" and switches towards "exasperated" exploration. A big question, which we have already partly addressed, is whether this analysis sits well with the Paretian expected utility theory or rather begets a different approach as proposed in (Peters, 2019).
- Do deviations from (Pareto) "rationality" and switches towards "exasperated" exploration depend on the dynamics of the optimization process as represented by Average Cumulative Reward observed over the limited number of trials available? Do deviations from Pareto correspond to a stressful situation? The observations of the values of ACR are compatible with this hypothesis.

In this paper we consider three measures $\sigma(x), h(x)$ or $z(x)$ to investigate which one would better fit with the experimental results in the framework of the Pareto rationality model (i.e., the one maximizing the number of choices laying on the associated Pareto frontier). A way to analyze this is a bi-objective optimization problem in the space of $\psi = (\zeta(x), u(x))$ where $\zeta(x)$ is the expected improvement and $u(x)$ represents the uncertainty. Research in Machine Learning has focused mostly on how to assess the informational utility of possible queries, we rather addressed the issue of the perception and quantification of probabilistic uncertainty itself. This problem is still an open question in Machine Learning and cognitive sciences and neither our results nor those prevailing in the rich literature about this issue provide unequivocal evidence about the algorithms underpinned by humans' decisions.

A very interesting perspective opened by these results stems from the observation that is that the locations in the space ($\psi = (\zeta(x), u(x))$) associated to Pareto-rational decisions and the degree of Pareto compliance do not change so much depending on the kernel, but rather on the uncertainty quantification: $\sigma(x), h(x)$ or $z(x)$. Moreover since $z(x)$ is most closely related to Pareto behaviour and it's also the one non kernel based this result could also be regarded as an indication that distance could be embedded directly in the learning paradigm.

We assume that two inputs $x$ and $x'$ are represented by two distributions $\delta$ and $\gamma$: many distributional distances can be used as Kullback-Leibler, Jensen-Shannon, $\chi$-square, total variation, Hellinger distance and Wasserstein which can be then "kernelized" as $K(x, x') = k(\delta, \gamma) = e^{-d(\delta, \gamma)}$.
But distance can also be used directly as the base of the learning process as first shown in (Balcan et al., 2008) Two remarks can be made:
- Using kernels limits the choice of distribution distances as the resulting kernel has to be definite positive: a widely used distance as Kullback-Leibler does not qualify. WST cannot be directly kernelized even if three recent papers (Le et al., 2019; Oh et al., 2019; De Plaen et al., 2020) have analysed the conditions in which WST yields a positive definite kernel.
- Working directly with distances could be computationally convenient dispensing with the computation of the kernel, its inversion, and the associated potential numerical instability.

A distance-based framework for learning with probability distributions has been proposed in (Rakotmamonjy et al., 2018). Based on the Wasserstein distance the paper contains a performance analysis of kernel versus distance-based classifiers.
A distributional distance-based learning has been shown to be very effective also in simulation-optimization problems over discrete structures (Ponti et al., 2021): the Multi Objective Evolutionary Algorithm based on the Wasserstein distance (MOEA/WST) has been shown to be more sample efficient than benchmark evolutionary approaches.

# 8. Declarations

### 8.1. Funding
We greatly acknowledge the DEMS Data Science Lab, Department of Economics Management and Statistics (DEMS), for supporting this work by providing computational resources.

### 8.2. Conflicts of interest/Competing interests (include appropriate disclosures)
Authors declare that they do not have any conflicts of interests or competing interests.

### 8.3. Availability of data and material (data transparency)
Both data and code for reproducing analysis and results of this paper are available at the following link: https://github.com/acandelieri/humans_strategies_analysis.

### 8.4. Code availability (software application or custom code)
Both data and code for reproducing analysis and results of this paper are available at the following link: https://github.com/acandelieri/humans_strategies_analysis.

### 8.5. Ethics approval (include appropriate approvals or waivers)
Informed consent was given in accordance with the university's procedure and the Helsinki declaration.

# Appendix A

## A1. Data collection

To collect data about humans' strategies we have used a gaming application based on the implementation used in (Candelieri et al., 2020). Figure 14 shows the web-based Graphical User Interface (GUI) of our game, with a game play example. The game field, with previous decisions and observations, as well as the score and remaining "shots", are reported. The game target is searching for the location having the highest score.

Fourteen volunteers have been enrolled (among colleagues and friends), asking for solving ten different tasks each (only for Game mode #1). Each task refers to a global optimization test function, which subjects "learn and optimize" by clicking at a location and observing the associated score (aka reward). For each task, every player has a maximum number of 20 clicks (decisions) available. The 10 global optimization test functions adopted are depicted in Appendix (A1). Since these functions are related to minimization problems, the score returned to the player is $-f(x)$, translating them into maximization tasks.

Finally, the game has been developed in R, specifically R-shiny for the web-based GUI. All the analytical components, described in the following, have been also developed in R as backend of the application.

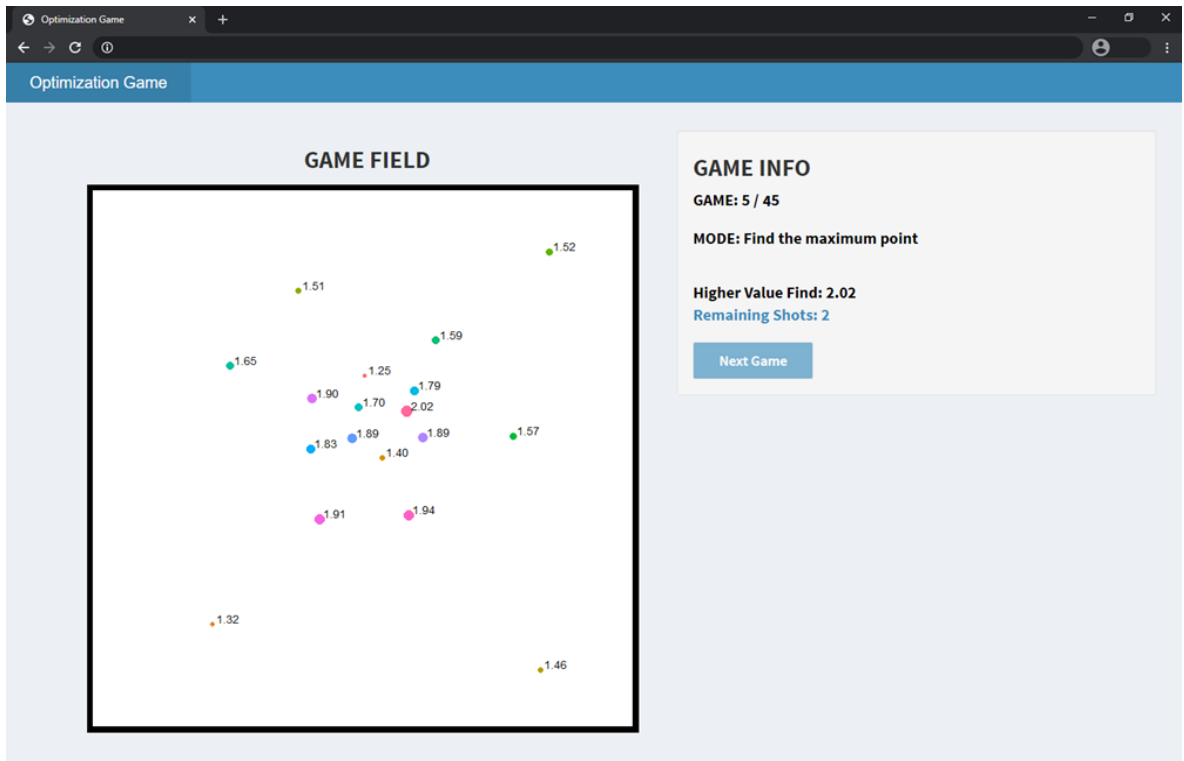

*Figure 14.* Web-based Graphical User Interface (GUI) of our game: a game play example

# Appendix B

## B1. The ten test problems

The ten global optimization test functions used in this study, including their analytical formulations, search spaces and information about optimums and optimizers, can be found at the following link:
https://www.sfu.ca/~ssurjano/optimization.html

Since they are minimization test functions, we have returned $-f(x)$ as score in order to translate them into the maximization problems depicted in Figure 15.

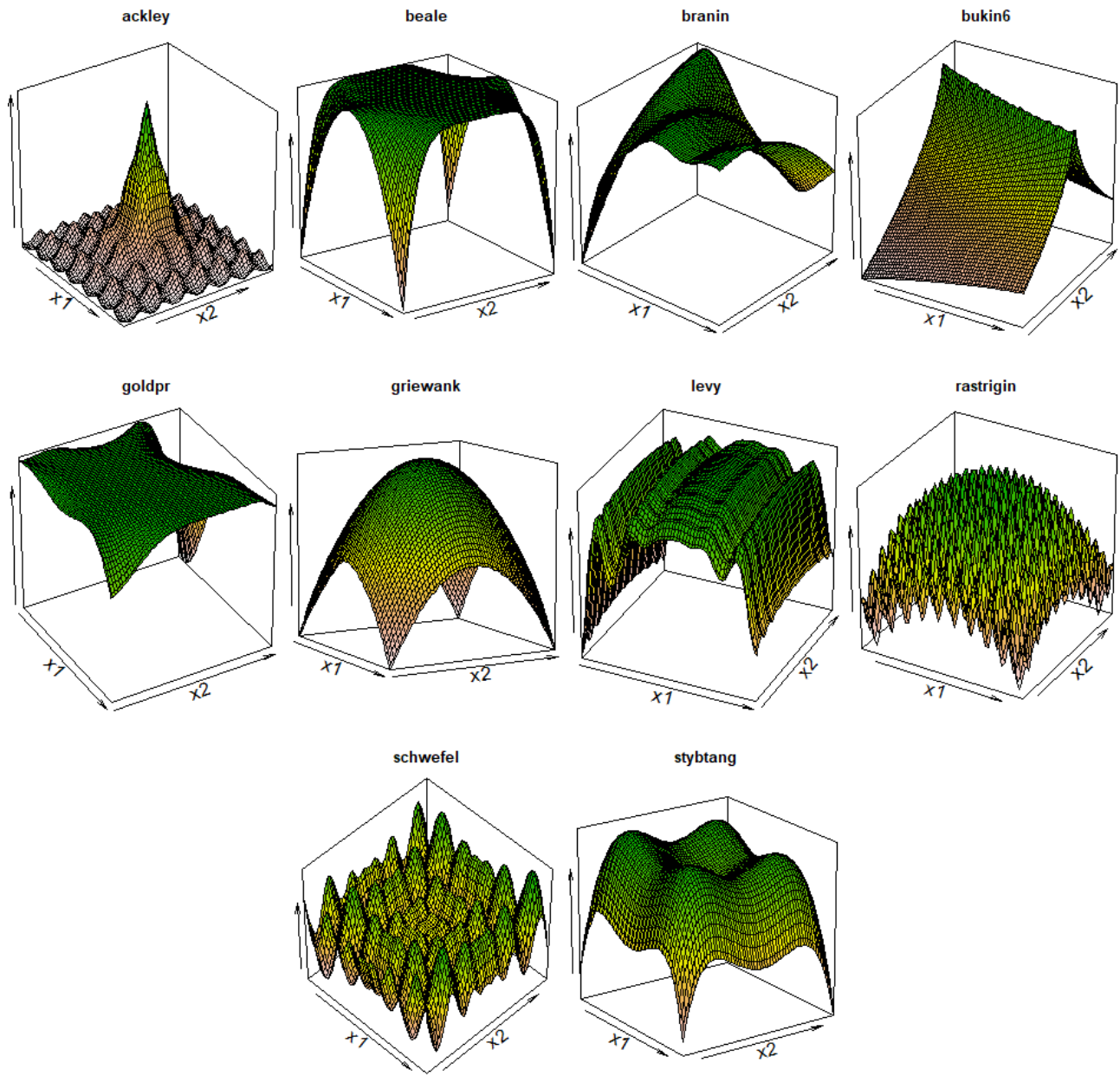

*Figure 15. The 10 test problems considered in this study.*

## B2. Distances from Pareto frontiers for each player, by test function

The following 10 figures – one for each test function – report the distances of each decision from the Pareto frontiers and for each player.

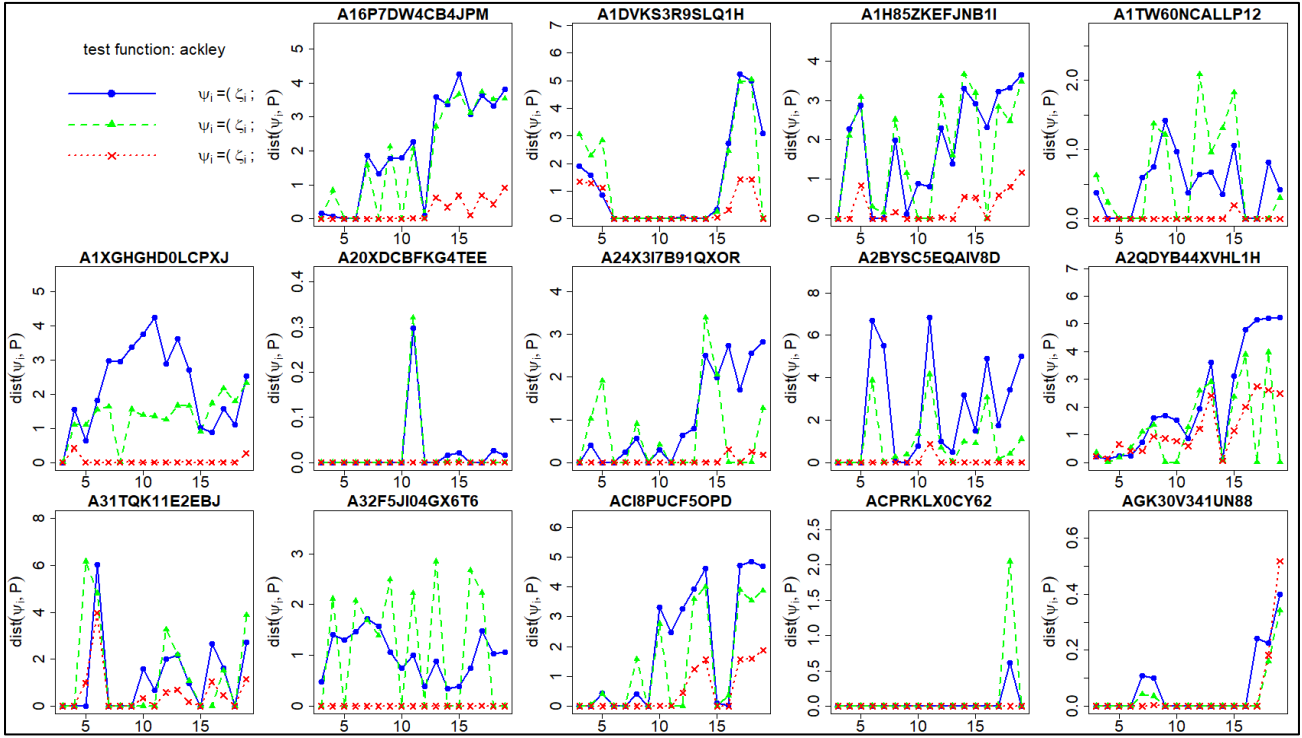

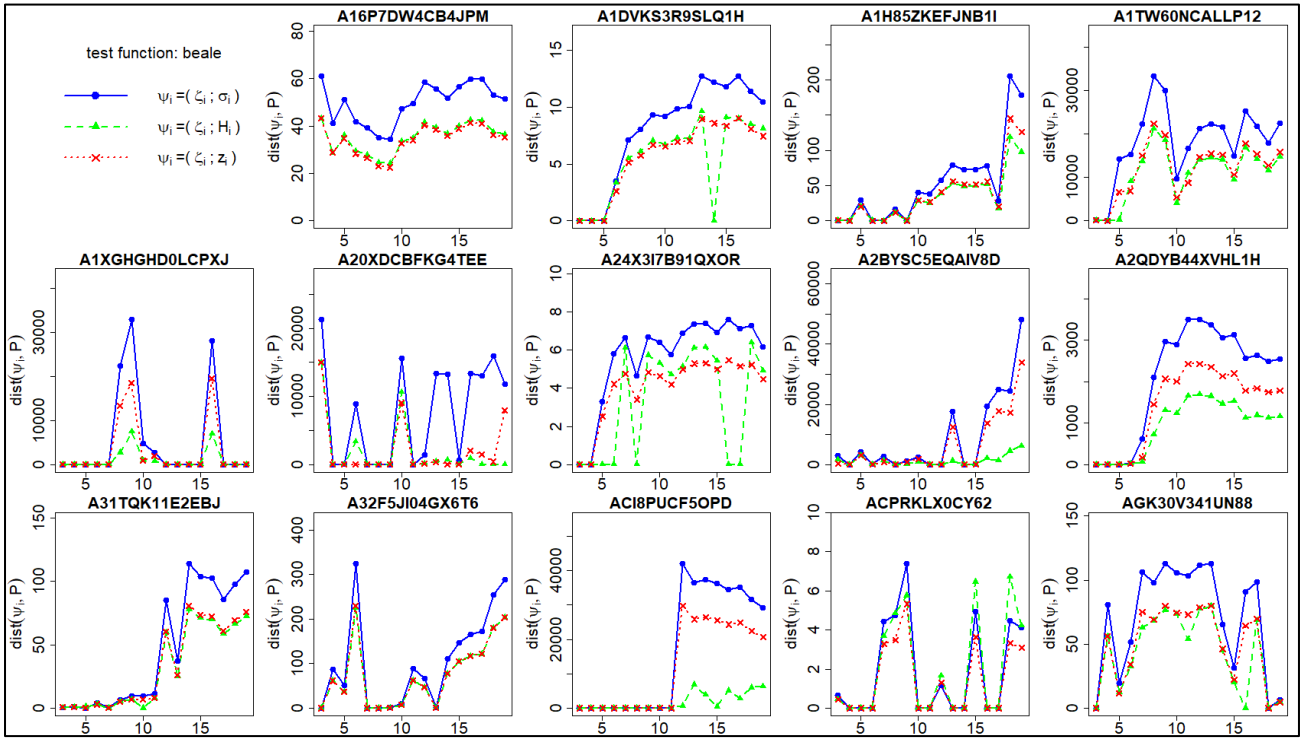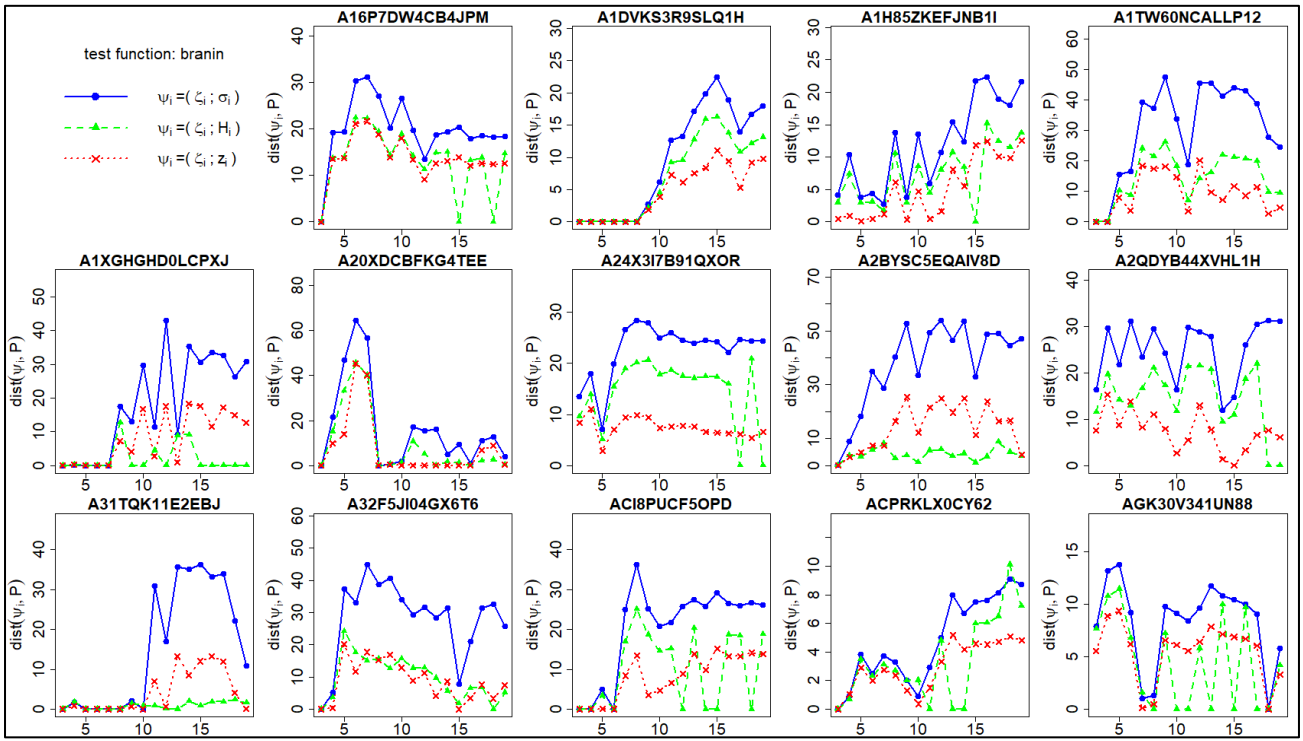

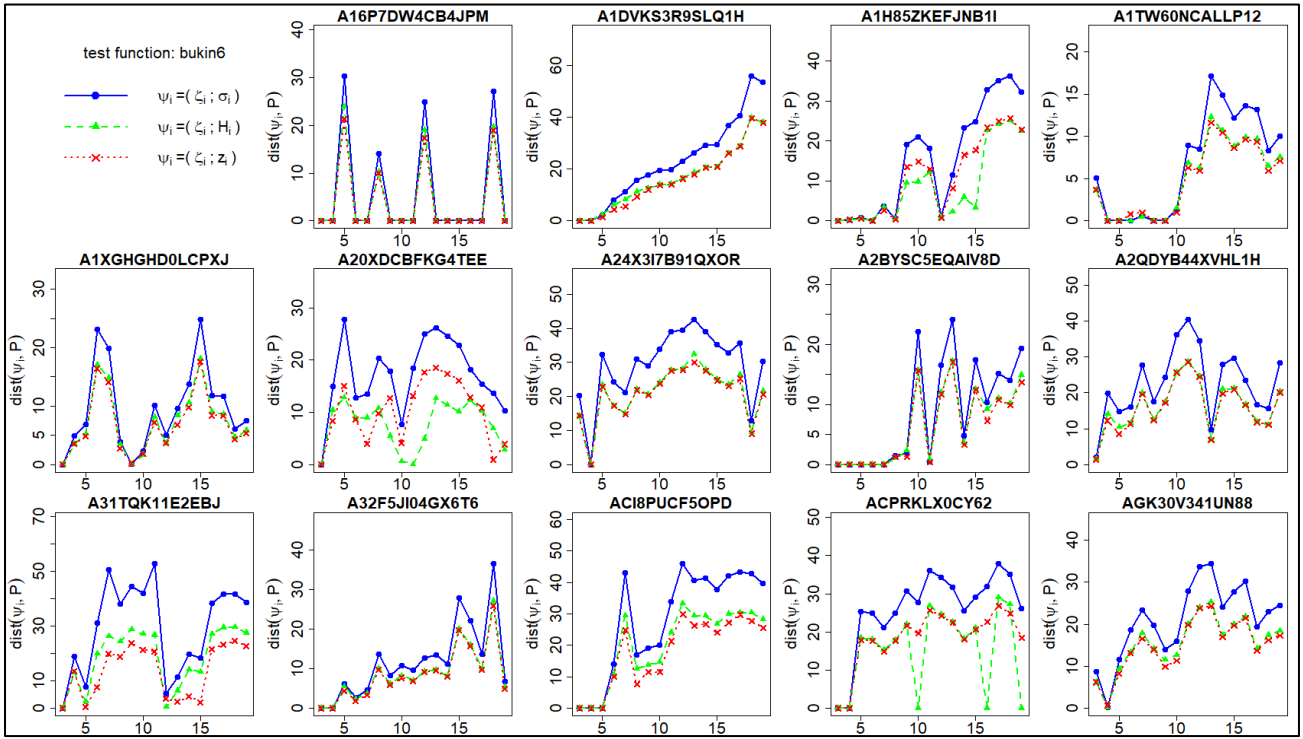

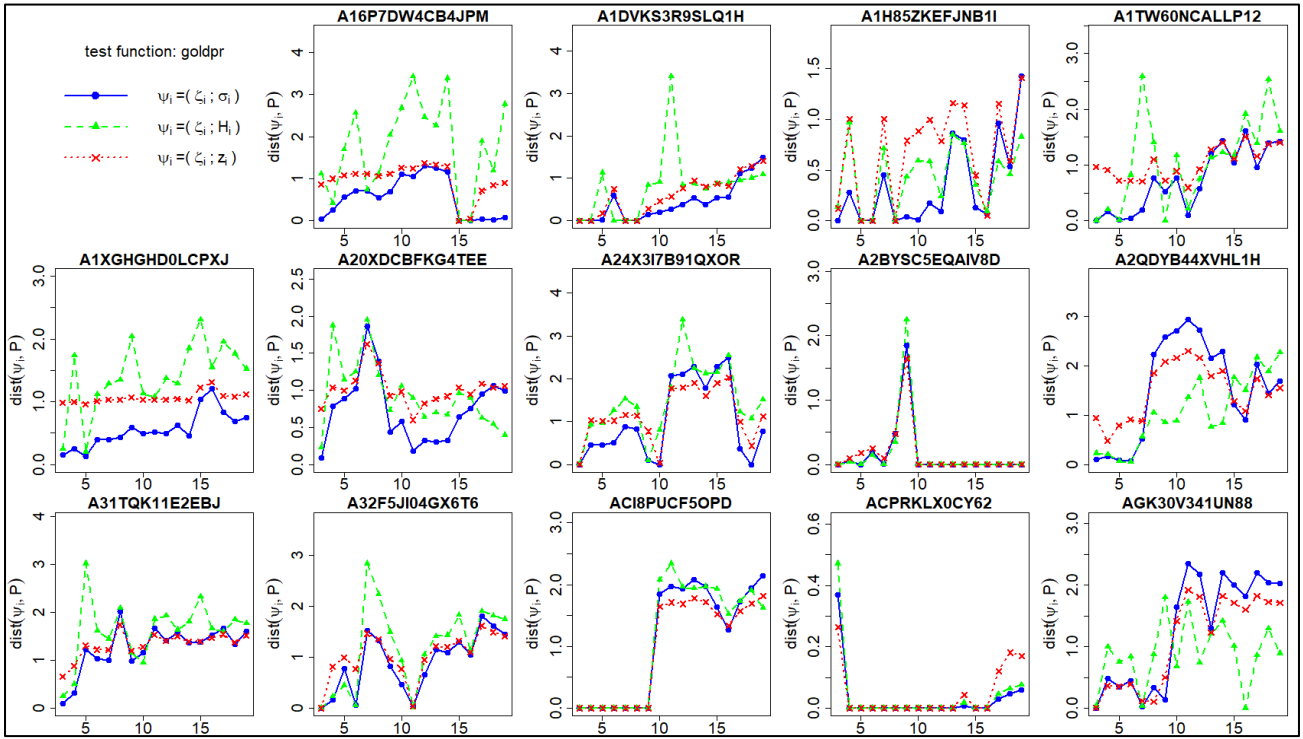

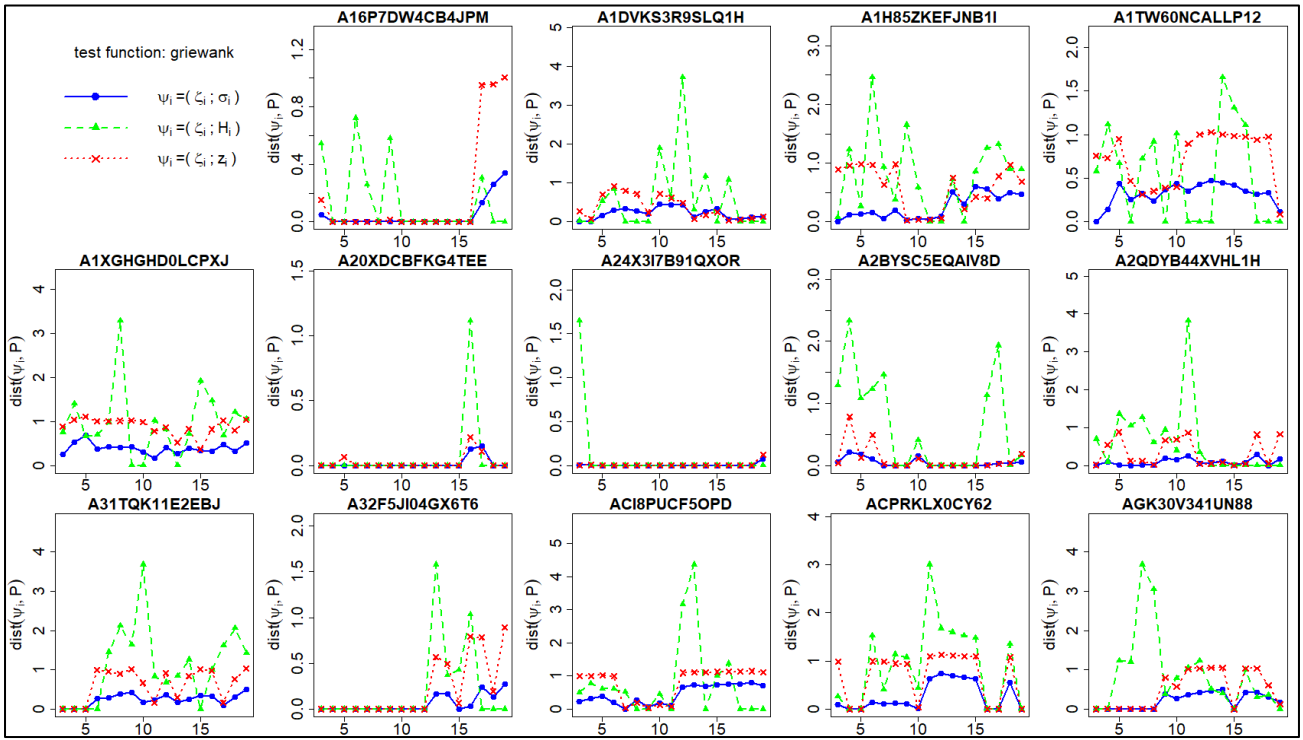

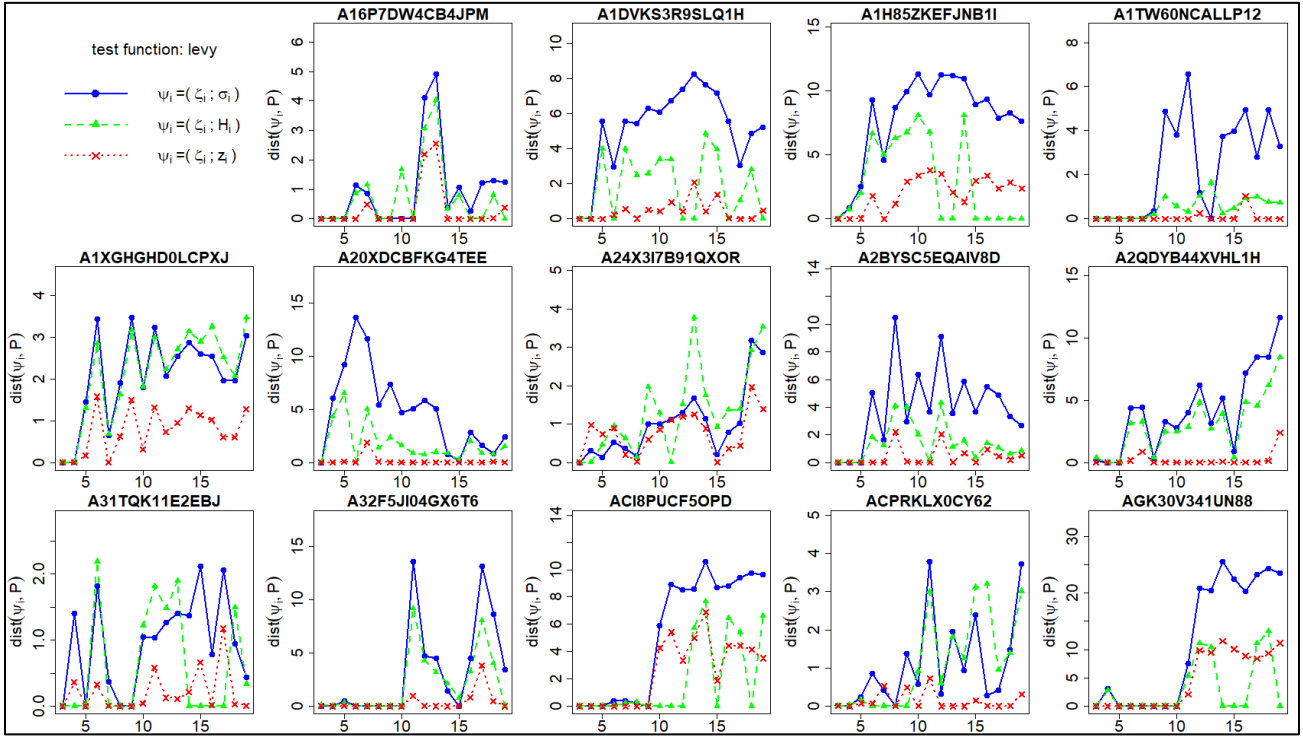

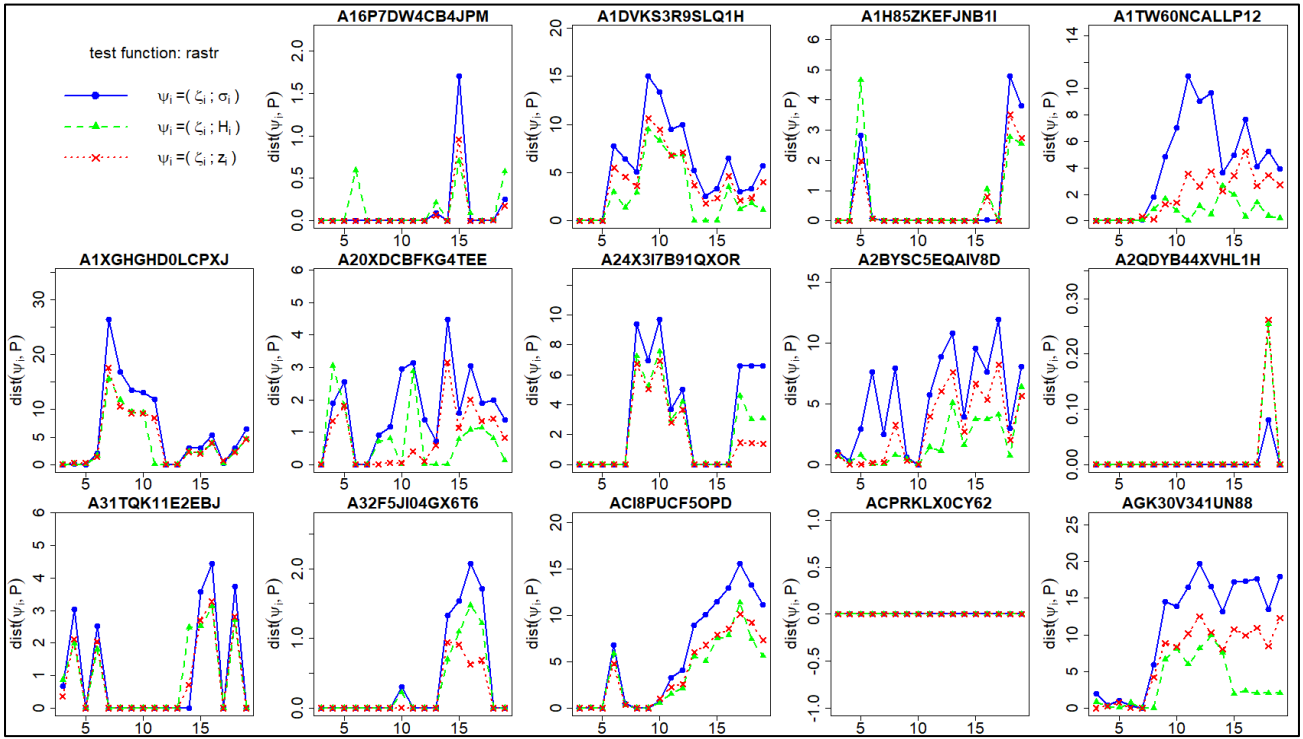

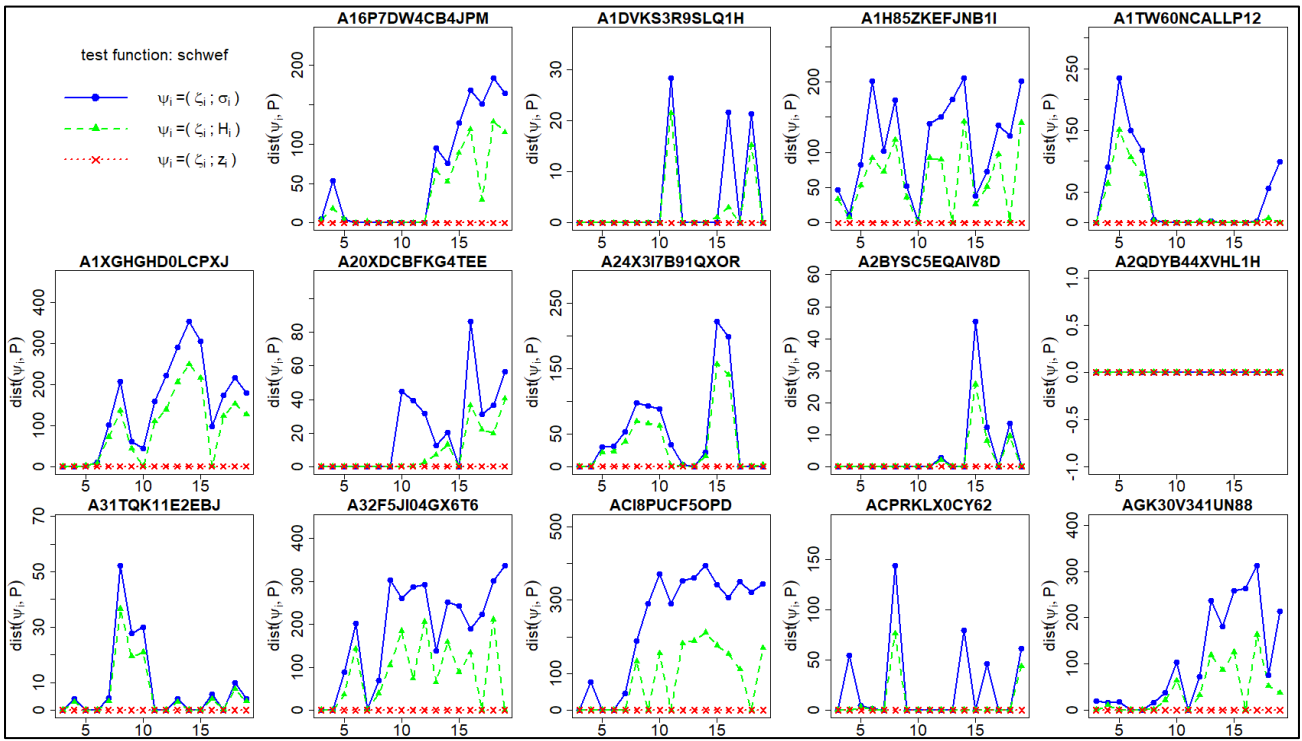
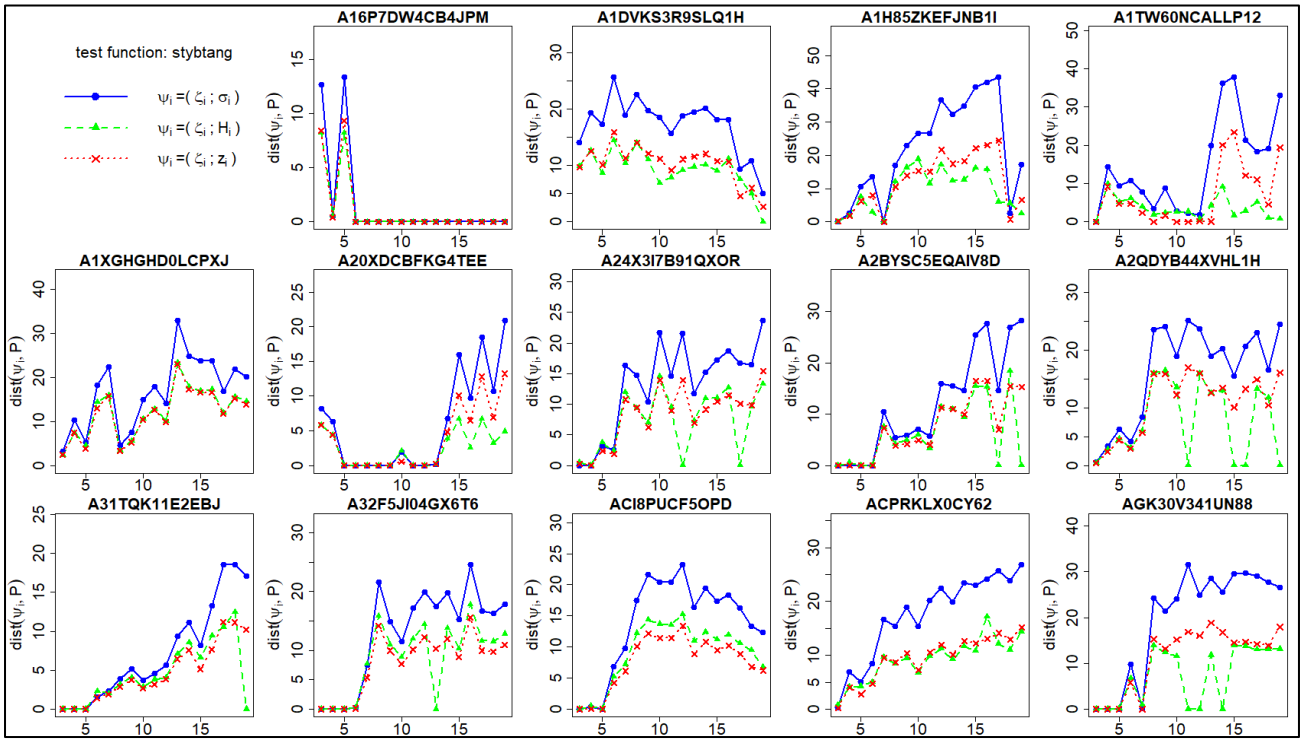

## B3. Distances from Pareto frontiers for each test functions, by player

The following 14 figures – one for each player – report the distances of each decision from the Pareto frontiers and with respect to each test function.

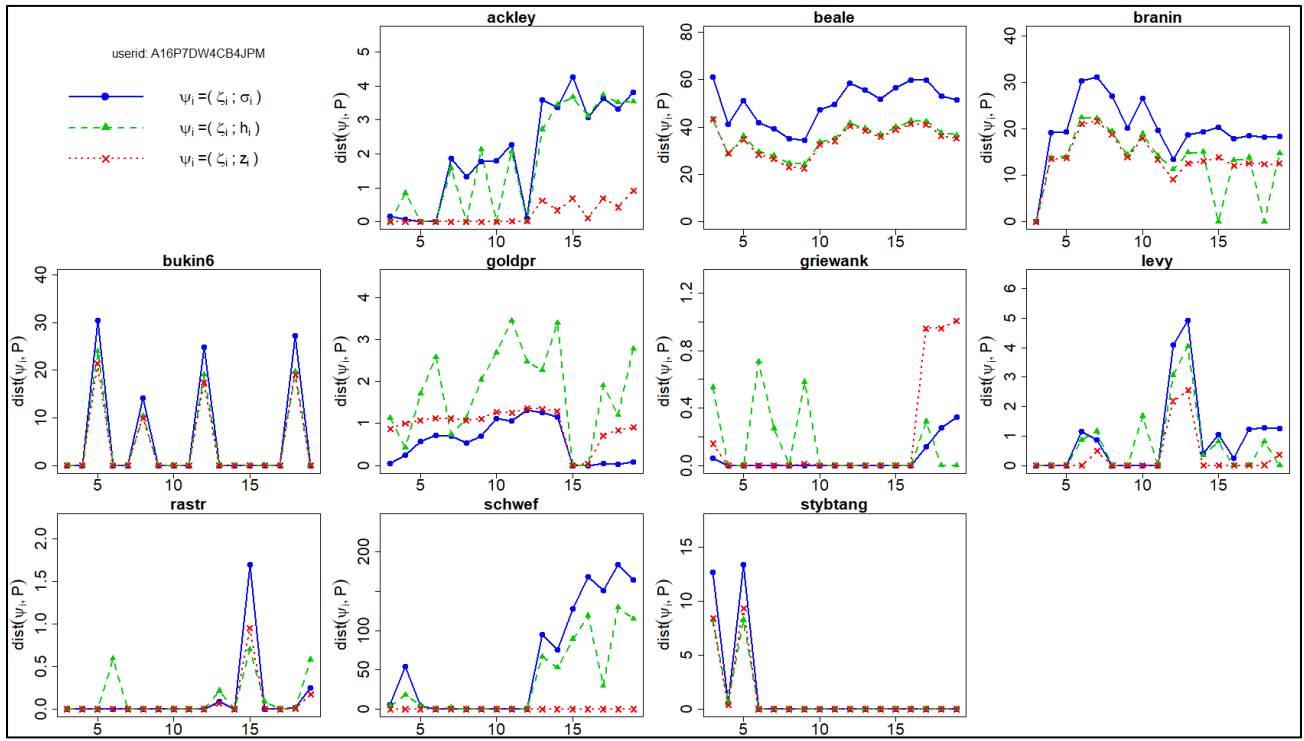

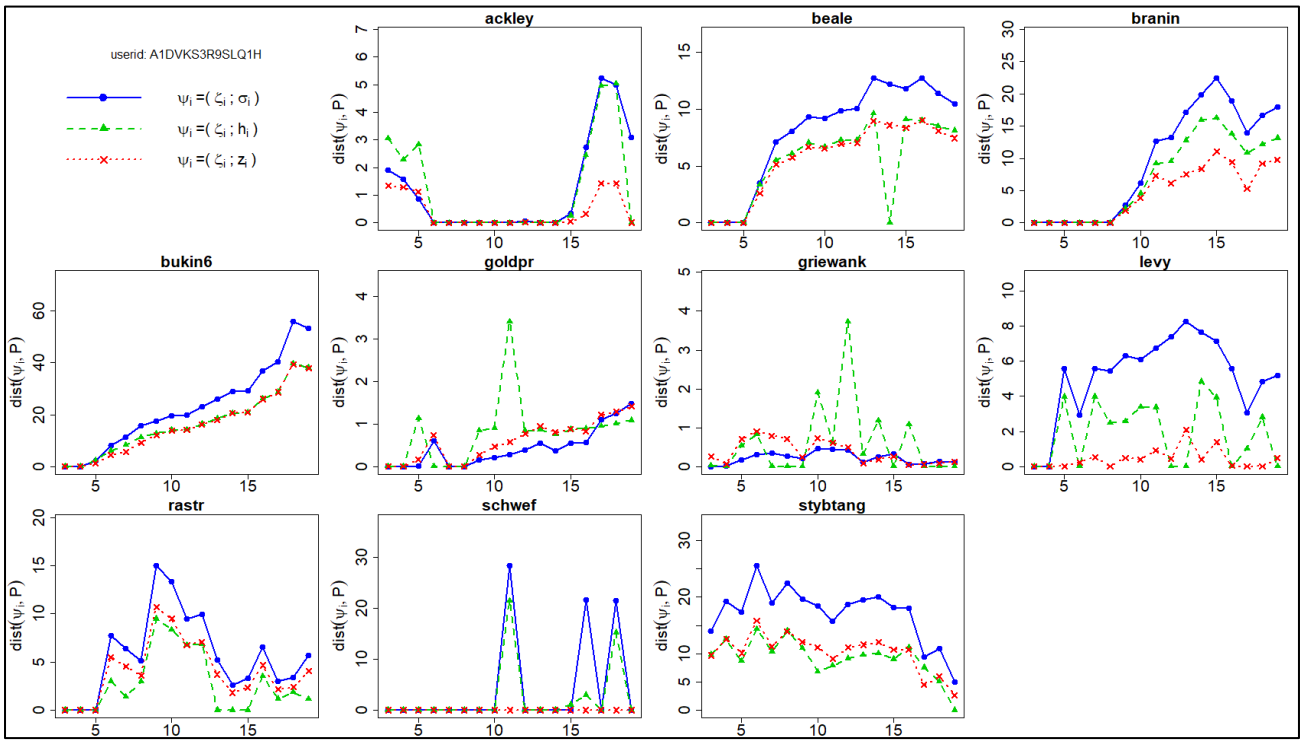

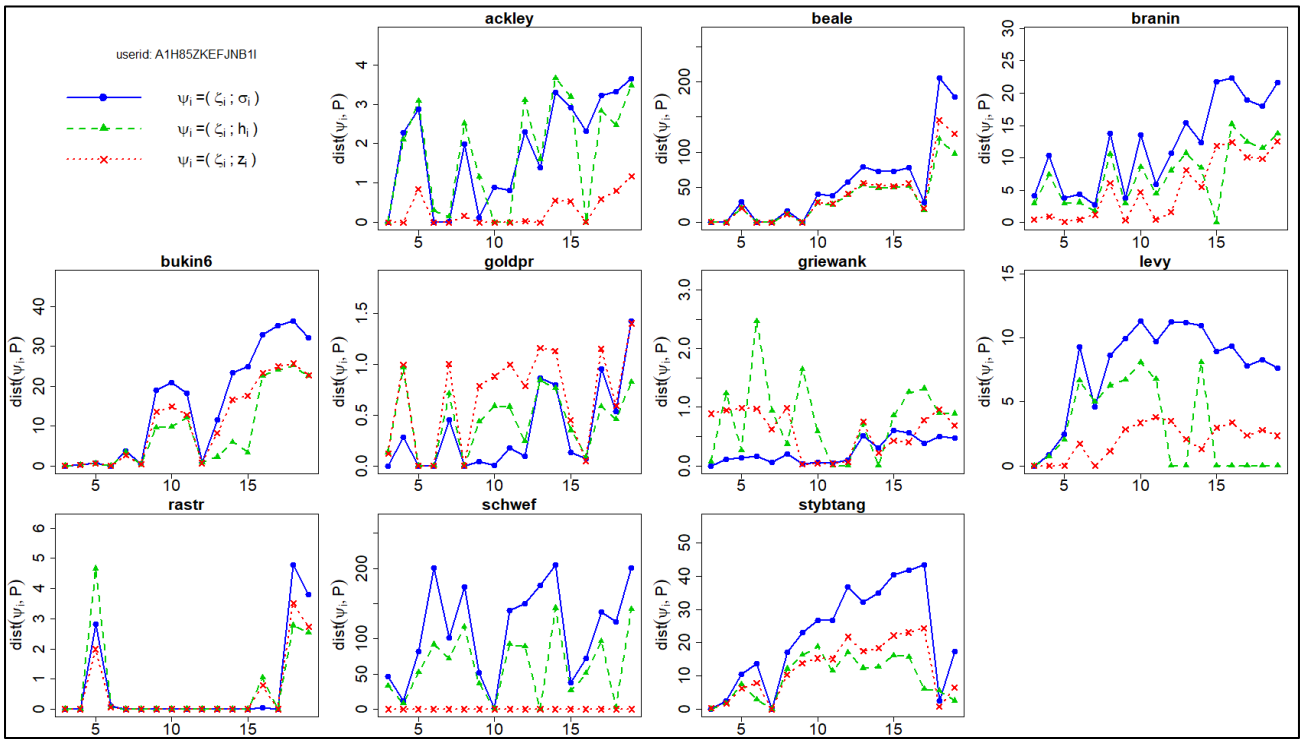

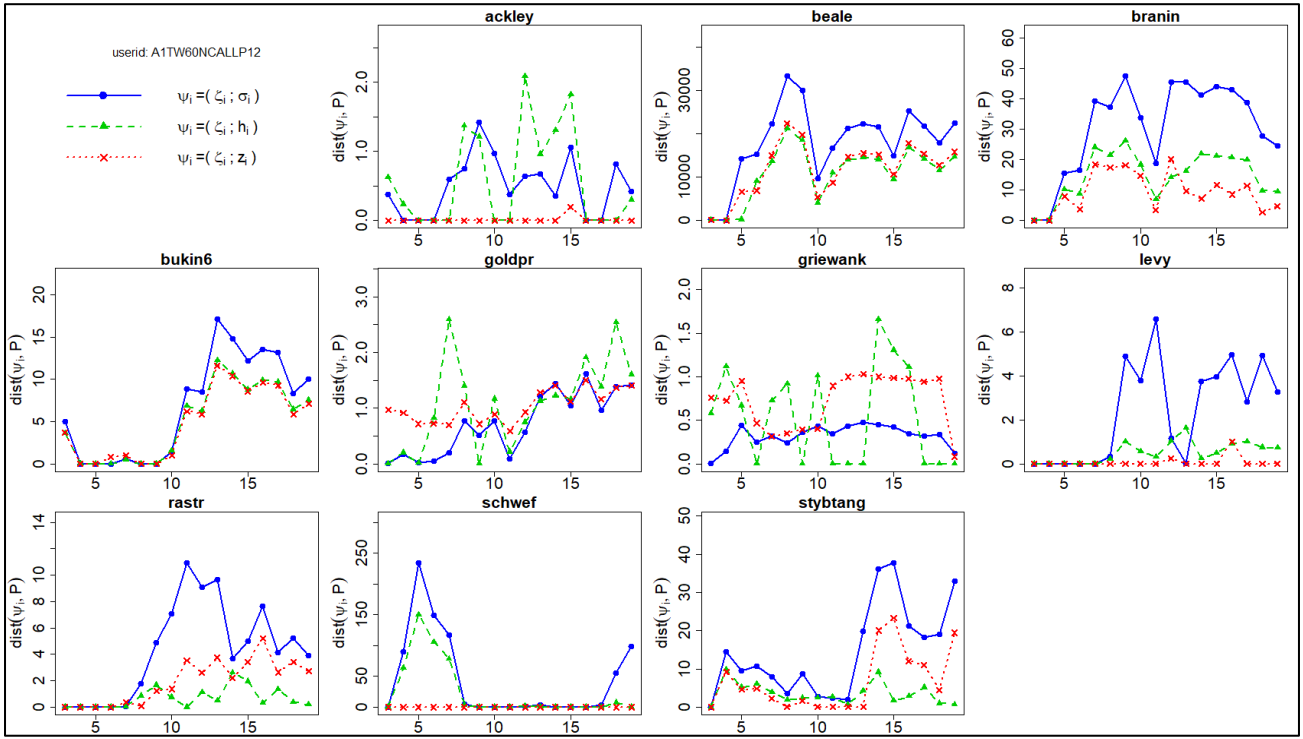

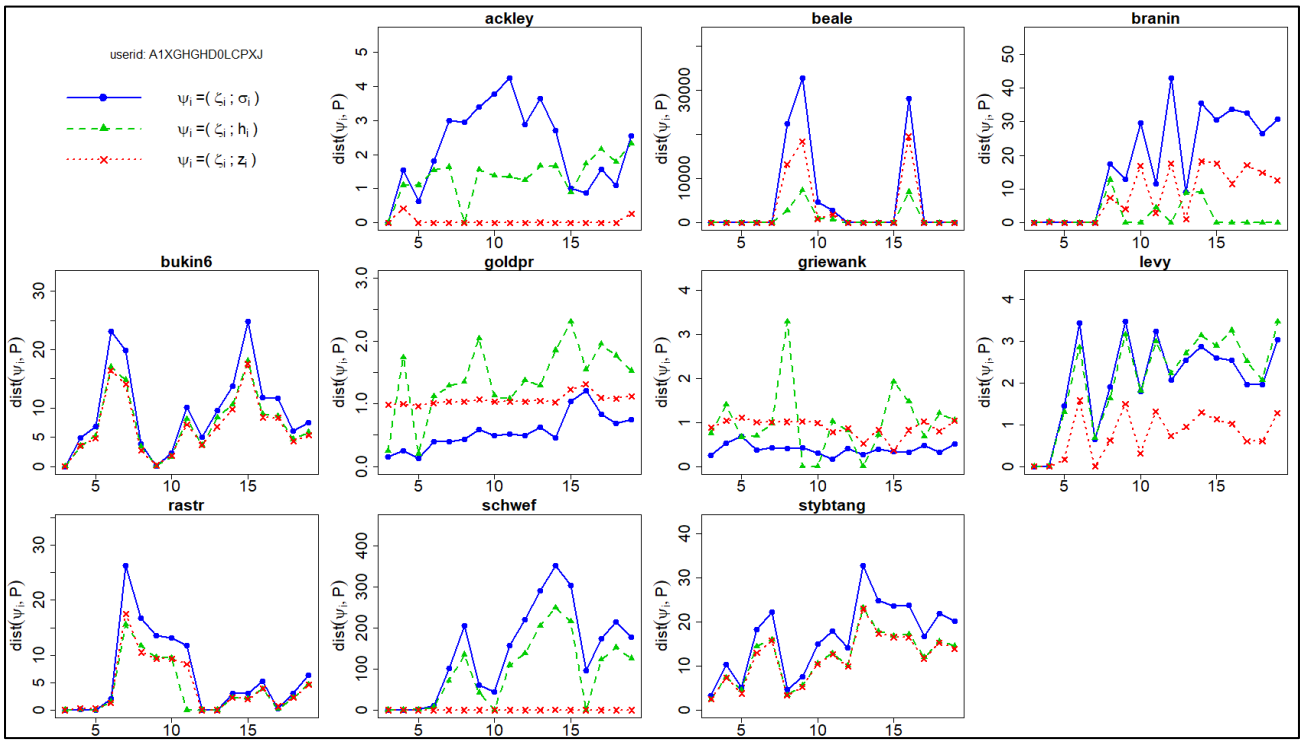
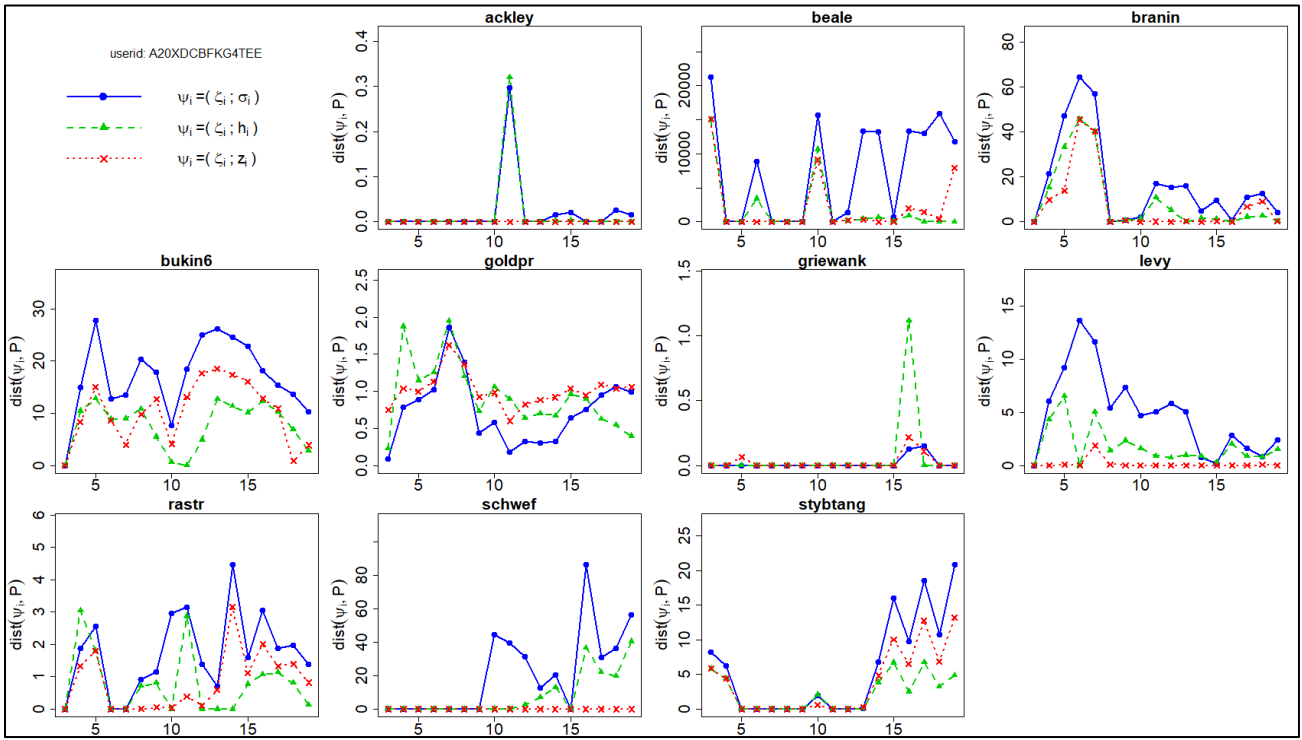

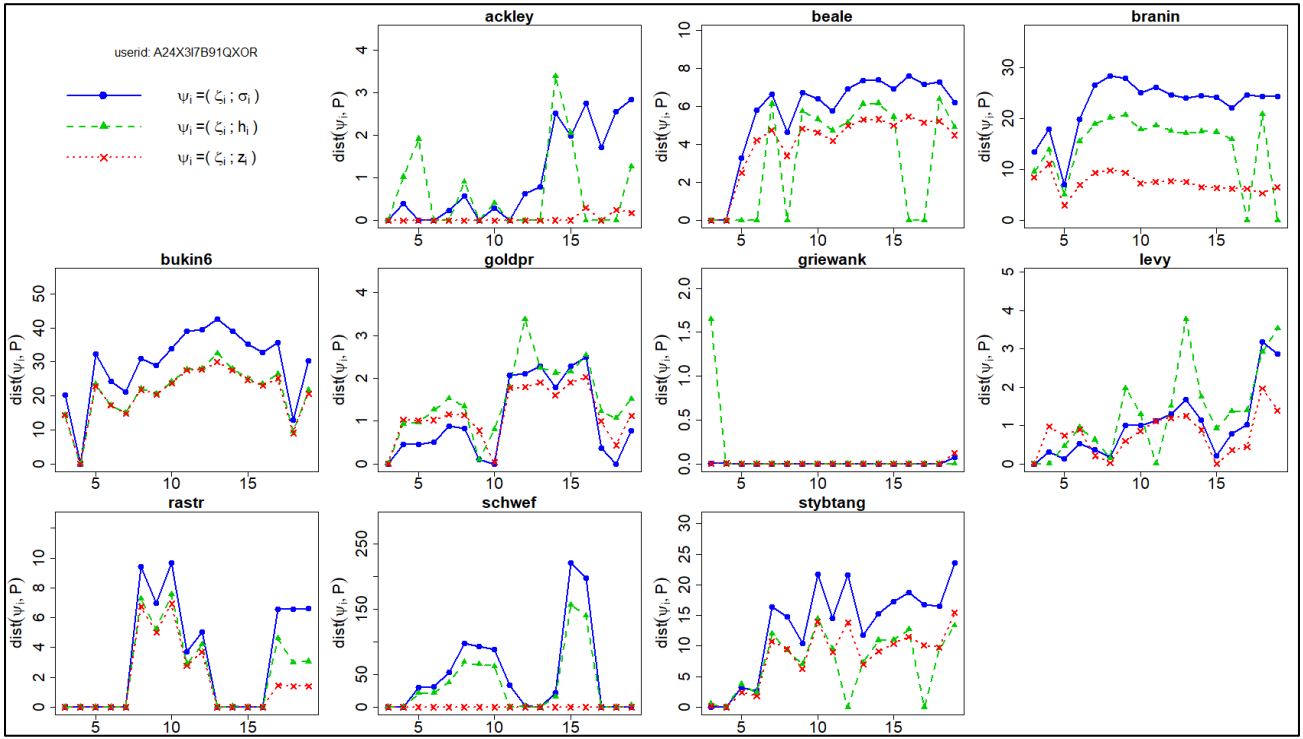

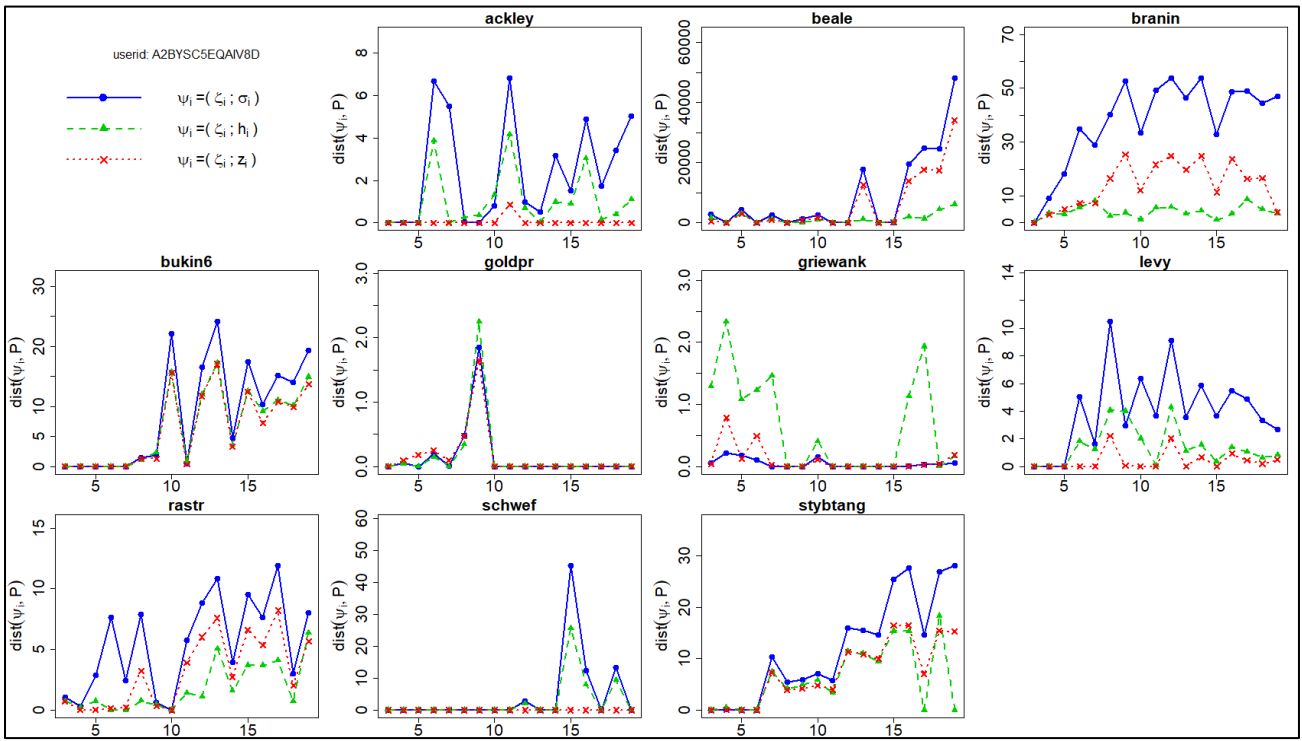

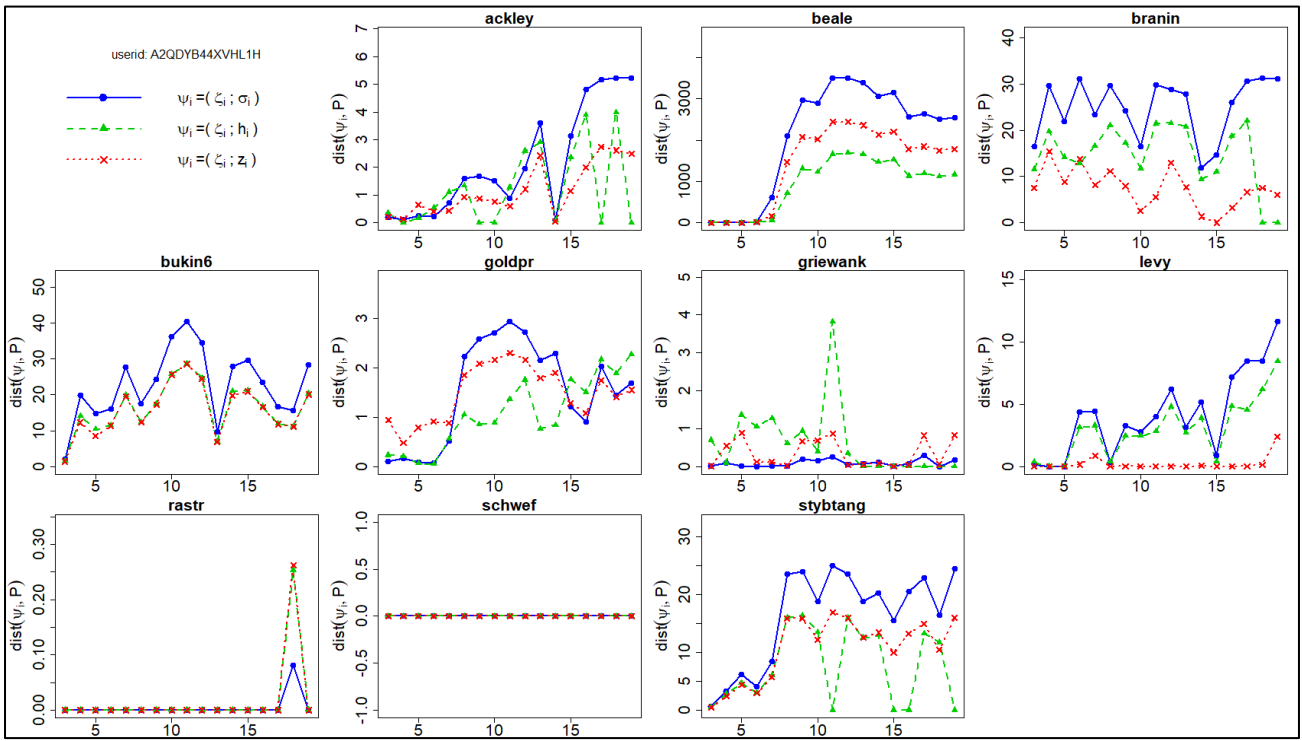

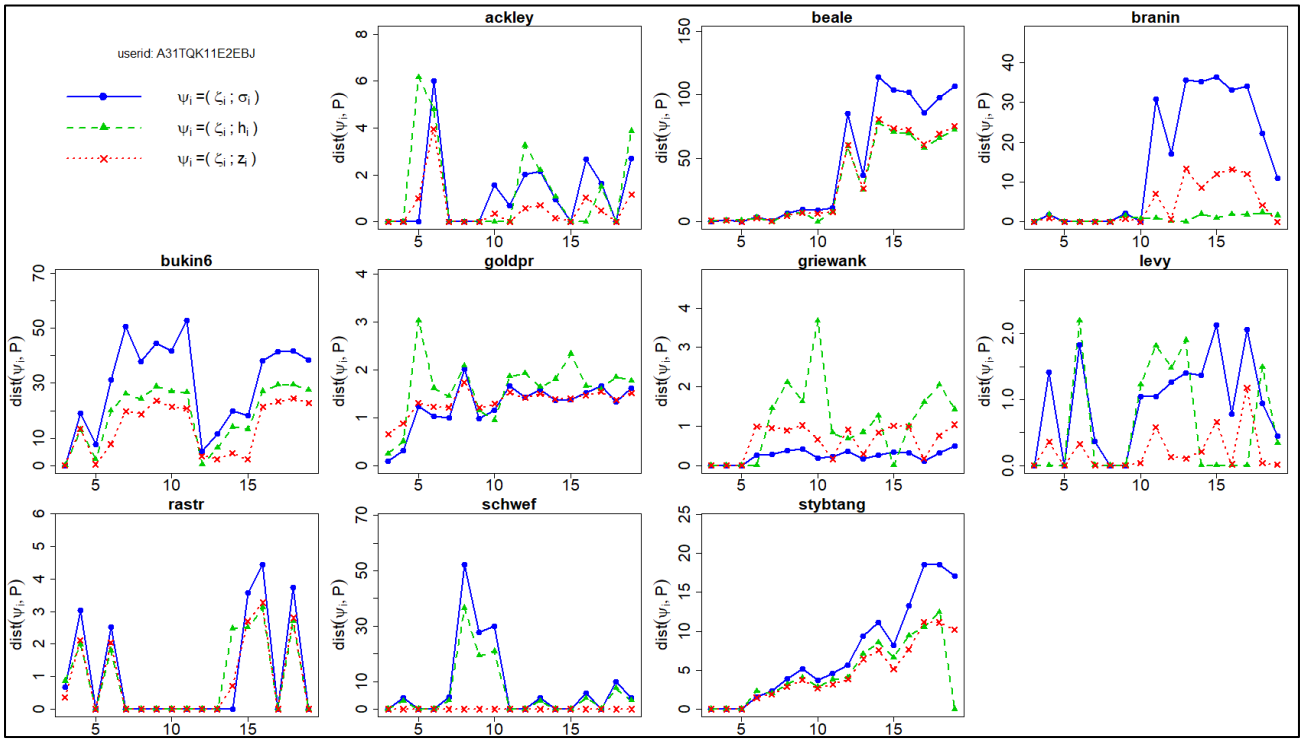

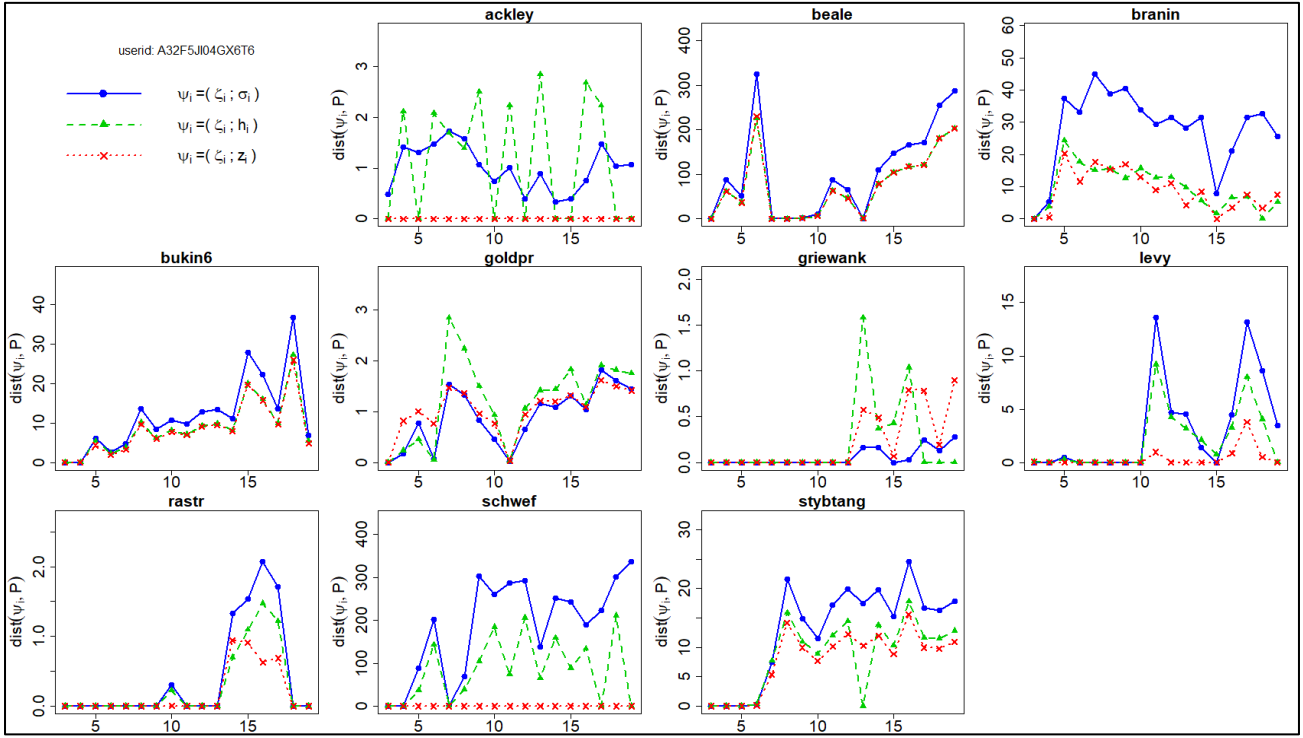

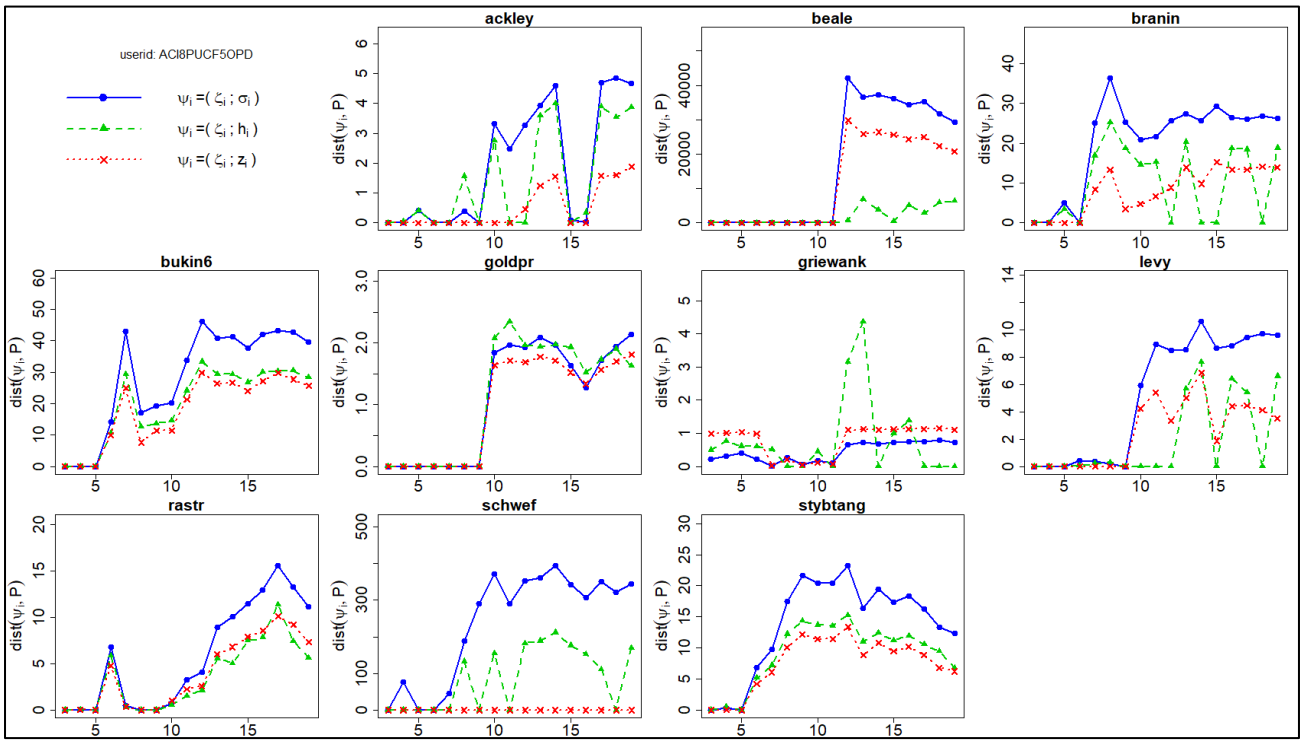

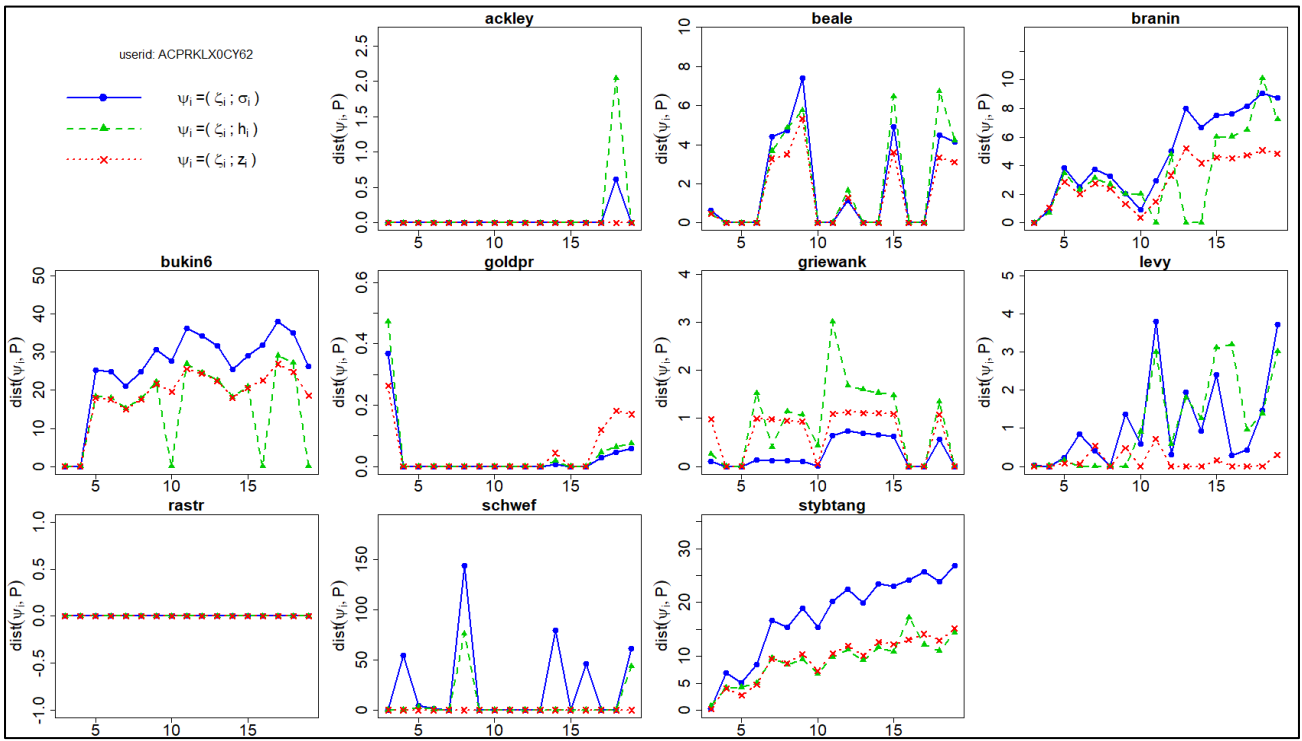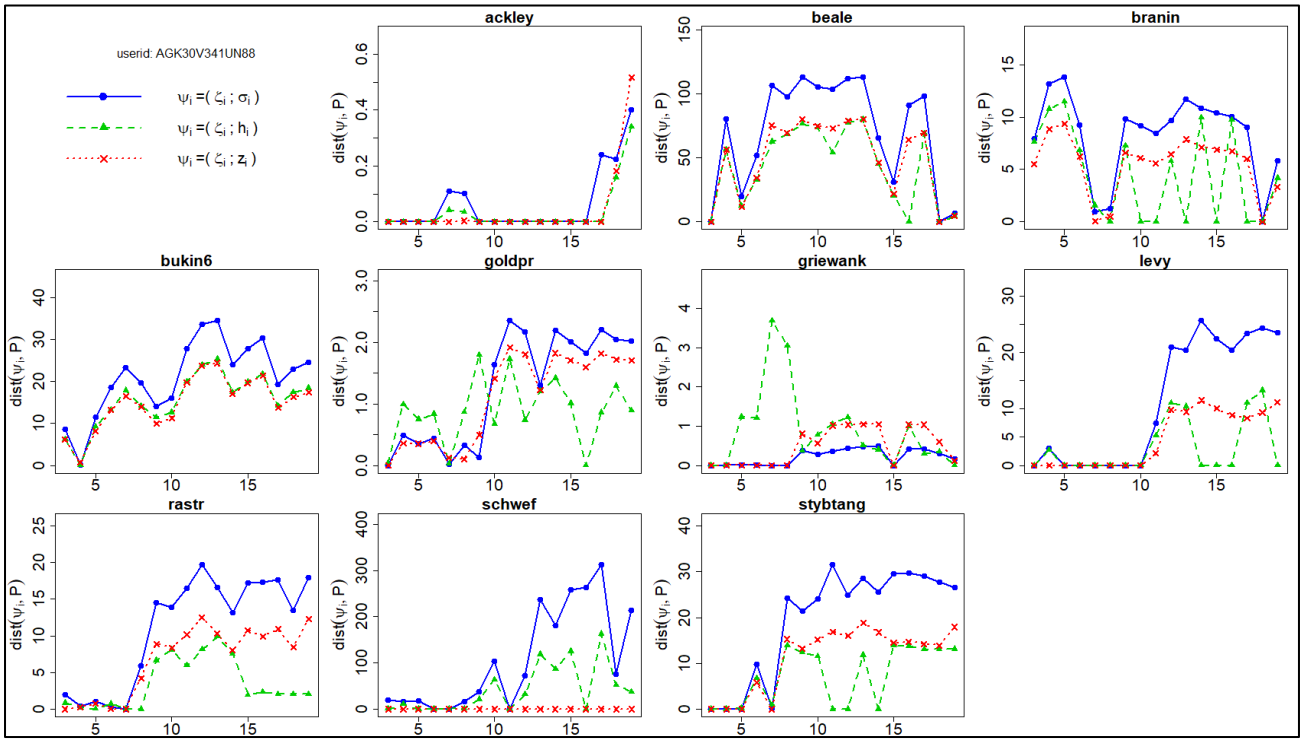